\documentclass[a4paper,12pt]{article}

\usepackage{amsmath,amssymb}

\allowdisplaybreaks

\begin{document}

\begin{center}

\textbf{\Large Fusion rule algebras related to a pair \\ 
of compact groups}

\medskip

\vspace{0,2 cm}

{\large Narufumi Nakagaki and Tatsuya Tsurii}

\medskip

{\large 15th  March 2017}

\end{center}

\begin{abstract}
The purpose of the present paper is to investigate a fusion rule algebra  
arising from irreducible characters of a compact group $G$ 
and a closed subgroup $G_0$ of $G$ with finite index.  
The convolution of this fusion rule algebra is introduced by inducing irreducible representations 
of $G_0$ to $G$ and by restricting irreducible representations of 
$G$ to $G_0$. 
%The method of proof relies on  
%character formulae of induced representations of compact  
%groups and of Frobenius' reciprocity theorem. 

\bigskip

\medskip



\end{abstract}

\bigskip

\section{Introduction}
Since an axiom of a hypergroup was established by  
C. F. Dunkl (\cite{D1}, \cite{D2}, 1973), 
R. I. Jewett (\cite{J}, 1975) and R. Spector  (\cite{S}, 1975), 
there are many researches about hypergroups. More precisely, refer to \cite{BH}. 
H. Heyer, S. Kawakami, T. Tsurii and S. Yamanaka introduced 
a hypergroup $\mathcal{K}(\widehat{G} \cup \widehat{G_0})$ arising 
from characters of a compact group $G$ and its closed subgroup $G_0$ 
with finite index (refer to [HKTY1]). 

\medskip

In the present paper, applying the results of [HKTY1], 
we consider a fusion rule algebra $\mathcal{F}(\widehat{G} \cup \widehat{G_0})$ 
arising from characters of a compact group and its closed subgroup with finite index.  
As in the case of [HKTY1], 
the method of constructions of fusion rule algebras $\mathcal{F}(\widehat{G} \cup \widehat{G_0})$ 
rely on the application of a character formula [H], 
Frobenius' reciprocity theorem [F] for compact groups, and 
character theory for induced representations of hypergroups [HKY].

\medskip

In Section 2 the preliminaries, we mention about an axiom of a 
countable discrete hypergroup and a fusion rule algebra. After the 
introduction, we consider relationships between hypergroups and fusion 
rule algebras, and introduce a notion of a fusion rule algebra join.

\medskip

In Section 3, we introduce a notion of a fusion rule algebra 
$\mathcal{F}(\widehat{G} \cup \widehat{G_0})$ arising from characters 
of a compact group $G$ and its closed subgroup $G_0$ with finite index. 
In this section, an admissible pair $(G, G_0)$  defined in [HKTY1] plays 
an important role to construct a fusion rule algebra 
$\mathcal{F}(\widehat{G} \cup \widehat{G_0})$. 

\medskip
In Section 4, we show some examples of fusion rule algebras 
$\mathcal{F}(\widehat{G} \cup \widehat{G_0})$. 
We note that hypergroups $\mathcal{K}(\widehat{G} \cup \widehat{G_0})$ 
in the sense of [HKTY1] is obtained by normalization of fusion rule 
algebras $\mathcal{F}(\widehat{G} \cup \widehat{G_0})$.

\section{Preliminaries} 

For a countable discrete set $K = \{c_0, c_1, c_2, \dots\}$, 
we denote the algebraic complex 
linear space based on $K$ together with finite support by $\mathbb{C}K$ , namely
\begin{align*}
\mathbb{C}K 
:= \left\{X = \sum_{j =0 }^\infty a_j {c_j}~:~
a_j \in \mathbb{C},~|\text{supp}(X)| < +\infty \right\},   
\end{align*}
where the support of $X$ is defined by 
$$
\text{supp}(X) := \{c_k : a_k \not = 0\}. 
$$

\medskip
\noindent
{\bf Axiom of a countable discrete hypergroup}   

\medskip

A {\it countable discrete hypergroup} $K = (K, \mathbb{C}K, \circ, *)$ consists 
of a countable discrete set $K = \{c_0, c_1, c_2, \dots\}$ 
together with an associative product (called convolution) $\circ$
and an involution $*$ in $\mathbb{C}K$ satisfying the following conditions. 

\begin{itemize}
\item[(H1)] The space ($\mathbb{C}K, \circ, *$)  
is an associative $*$-algebra with unit ${c_0}$. 
\item[(H2)] For $c_i, c_j \in K$, the convolution ${c_i} \circ {c_j}$ 
belongs to $\mathbb{C}K$, the coefficients of ${c_i} \circ {c_j}$ are 
non-negative real numbers and sum of them is one. 
\item[(H3)] $K^* = K$ i.e. $c_i^* \in K$ for $c_i \in K$. 
Moreover $c_j = c_i^*$ if and only if 
$c_0 \in \text{supp}({c_i} \circ {c_j})$. 
\end{itemize}
A countable discrete hypergroup $K$ is said to be 
commutative if the convolution $\circ$ is commutative. 

\medskip

By the above axiom (H3), 
\begin{align*}
c_i \circ c_i^* = a_0 c_0 + \cdots~~~(a_0 \not = 0). 
\end{align*}
Then the weight $w(c_i)$ of $c_i$ is defined by 
$$
w(c_i) := \frac{1}{a_0}. 
$$ 

\medskip
\noindent
{\bf Axiom of a fusion rule algebra}

\medskip

A {\it fusion rule algebra} $F= (F, \mathbb{C}F, \circ, *)$ 
consists of a countable discrete set $F =\{X_0, X_1, X_2, \dots\}$ 
together with an associative product (called convolution) $\circ$
and an involution $*$ in $\mathbb{C}F$ satisfying the following conditions. 
\begin{itemize}
\item[(F1)] The space ($\mathbb{C}F, \circ, *$)  
is an associative $*$-algebra with unit ${X_0}$.
\item[(F2)] For $X_i, X_j \in F$, the convolution ${X_i} \circ {X_j}$ 
belongs to $\mathbb{C}F$ and the coefficients of ${X_i} \circ {X_j}$ are 
non-negative integers. 
\item[(F3)] $F^* = F$ i.e. $X_i^* \in F$ for $X_i \in F$. 
Moreover $X_j = X_i^*$ if and only if 
$X_0 \in \text{supp}({X_i} \circ {X_j})$ and 
$$
{X_i} \circ {X_j} = X_0 + \sum_{\substack{X_k \in S(X_i, X_j) \\ X_k \not = 0}}  a_{ij}^k X_k. 
$$
where $S(X_i, X_j) := \text{supp}(X_i \circ X_j)$.  
\end{itemize}
A fusion rule algebra $F$ is said to be commutative if the convolution $\circ$ is commutative.

\medskip
\noindent
{\bf A relationship between  hypergroups and fusion rule algebras}

\medskip

Let $F$ be a fusion rule algebra. The mapping $d$ from $F$ to $\mathbb{R}_+^\times = (0, \infty)$ 
is called a {\it dimension function} of $F$ if $d$ is a homomorphism in the sense that 
\begin{align*}
d(X_i) d(X_j) = \sum_{X_k \in S(X_i,X_j)} a_{ij}^k d(X_k)~~
\text{where}~~
X_i \circ X_j = \sum_{X_k \in S(X_i,X_j)} a_{ij}^k X_k. 
\end{align*}
Then $d$ is uniquely extendable as a linear mapping from $\mathbb{C}F$ to $\mathbb{C}$ 
and satisfies 
$$
d(X_i \circ X_j) = d(X_i) d(X_j).
$$
For $X_i \in F$, put 
$$
c_i := \frac{1}{d(X_i)}X_i. 
$$
Then the set $K = \{c_0, c_1, c_2,\dots\}$ becomes a countable discrete hypergroup. 
We call this hypergroup $K$ a hypergroup which is normalized $F$, 
we write $\mathcal{K}_d(F)$.

\medskip
Let $F = \{X_0, X_1, \dots, X_n\}$ be a finite fusion rule algebra 
and $\mathcal{K}_d(F) = \{c_0, c_1, \dots \\, c_n\}$ a hypergroup which normalized $F$. 
The convolution of $c_i \circ c_i^*$ is 
\begin{align*}
c_i \circ c_i^* = \frac{1}{d(X_i)d(X_i^*)}X_i \circ X_i^* 
= \frac{1}{d(X_i)^2}X_0 + \cdots 
= \frac{1}{d(X_i)^2}c_0 + \cdots.  
\end{align*}
Hence the weight of $c_i$ is 
\begin{align*}
w(c_i) = d(X_i)^2. 
\end{align*}
Put 
$$
R(F) = \sum_{k=0}^n d(X_k) X_k. 
$$
Then the following proposition holds. 

\medskip
\noindent
{\bf Proposition 2.1 } For $X_i \in F$, 
\begin{align*}
X_i \circ R(F) = d(X_i) R(F). 
\end{align*}

\medskip
\noindent
{\bf Proof } Since $F= \{X_0, X_1, \dots, X_n\}$ is a fusion rule algebra, 
there exists a hypergroup $\mathcal{K}_d(F) = \{c_0, c_1, \dots,c_n\}$ 
which normalized $F$, the elements $c_i$ of $\mathcal{K}_d(F)$ are written by
$$
c_i = \frac{1}{d(X_i)}X_i. 
$$ 
Hence, 
\begin{align*}
R(F) = \sum_{k=0}^n d(X_k)^2 c_k = \sum_{k=0}^n w(c_k) c_k. 
\end{align*}
Therefore, $R(F)$ is the Haar measure of $\mathcal{K}_d(F)$, then
\begin{align*}
X_i \circ R(F) = d(X_i) c_i \circ R(F) = d(X_i) R(F).
\end{align*} 
\hspace{\fill}$\Box$

\bigskip
\noindent
{\bf Fusion rule algebra join}

\medskip
Let $H=\{X_0, X_1, \dots, X_n\}$ be a finite fusion rule algebra 
and $L = \{Y_0, Y_1\}$ a cyclic group of order two. The convolution and involution of 
$H$ and $L$ are $\circ_H$, $*_H$ and $\circ_L$, $*_L$ respectively. 
We assume that for all $X_i \in H$ 
the values of dimension function $d(X_i)$ are natural number.      
For the set 
$$
H \vee_F L = \{X_0, X_1, \dots, X_n, Y_1\}
$$ 
we define the convolution $\circ_F$ on $\mathbb{C}(H \vee_F L)$ by 
\begin{align*}
&X_i \circ_F X_j := X_i \circ_H X_j,\\
&X_i \circ_F Y_1 := d(X_i) Y_1,\\
& Y_1 \circ_F Y_1 :=  R(H).  
\end{align*}  
The involution $*_F$ defined by 
\begin{align*}
X_i^{*_F} := X_i^{*_H}~~{\rm and}~~Y_1^{*_F} := Y_1^{*_L} = Y_1. 
\end{align*}

\medskip
\noindent
{\bf Proposition 2.2 } $H \vee_F L = (H \vee_F L,\mathbb{C}(H \vee_F L), \circ_F, *_F)$ 
becomes a fusion rule algebra. 

\medskip
\noindent
{\bf Proof } By the definition of the convolution $\circ_F$ and the assumption
of a fusion rule algebra join,  
it is sufficient to show $H \vee_F L$ satisfies following associativity relations. 
\begin{itemize}
\item[(1)] $(X_i \circ_F X_j) \circ_F X_k = X_i \circ_F (X_j \circ_F X_k)$, 
\item[(2)] $(X_i \circ_F X_j) \circ_F Y_1 = X_i \circ_F (X_j \circ_F Y_1)$, 
\item[(3)] $(X_i \circ_F Y_1) \circ_F Y_1 = X_i \circ_F (Y_1 \circ_F Y_1)$, 
\item[(4)] $(Y_1 \circ_F Y_1) \circ_F Y_1 = Y_1 \circ_F (Y_1 \circ_F Y_1)$. 
\end{itemize}

(1) is clear because $H$ is a fusion rule algebra. 

\medskip

(2) 
\begin{align*}
(X_i \circ_F X_j) \circ_F Y_1 = \sum_{s=0}^n a_{ij}^s X_s \circ_F Y_1 
= \sum_{s=0}^n a_{ij}^s d(X_s) Y_1. 
\end{align*}
On the other hand, 
\begin{align*}
X_i \circ_F (X_j \circ_F Y_1) &= X_i \circ_F d(X_j) Y_1 = d(X_i)d(X_j) Y_1\\
&= d(X_i \circ_F X_j) Y_1 = d \left( \sum_{s=0}^n a_{ij}^s X_s \right) Y_1
= \sum_{s=0}^n a_{ij}^s d(X_s) Y_1. 
\end{align*}
Hence, the associativity relation (2) holds. 

\medskip
(3)
\begin{align*}
(X_i \circ_F Y_1) \circ_F Y_1 = d(X_i) Y_1 \circ_F Y_1 = d(X_i) R(H). 
\end{align*}
On the other hand,
\begin{align*}
X_i \circ_F (Y_1 \circ_F Y_1) = X_i \circ_F R(H). 
\end{align*}
By Proposition 2.1, the associativity relation (3) holds. 

\medskip
(4) is easy to check by simple calculation. 
Hence $H \vee_F L$ becomes a fusion rule algebra. 

\hspace{\fill}$\Box$

We call this algebra $H \vee_F L$ a fusion rule algebra {\it join} of $H$ by $L$.

\section{Fusion rule algebras related to admissible pairs}

\medskip

Let $G$ be a compact group which satisfies the second axiom of 
countability and let $\widehat{G}$ be the set of all equivalence classes of 
irreducible representations of $G$. Then $\widehat{G}$ is (at most  countable)  
discrete space which we write explicitly as  
$$
\widehat{G} = \{\pi_0, \pi_1, \dots, \pi_n,\dots \}, 
$$
where $\pi_0$ is the trivial representation of $G$. 
We denote by ${\rm Rep}^{\rm f}(G)$ the set of equivalence classes of 
finite-dimensional representations of $G$. For  $\pi \in {\rm Rep}^{\rm f}(G)$  
we consider the character of $\pi$ given by 
$$
Ch(\pi)(g) = {\rm tr}(\pi(g)),  
$$
for all $g \in G$.   
Put 
$$
\mathcal{F}(\widehat{G}) = \{Ch(\pi) : \pi \in \widehat{G}\}. 
$$
Then $\mathcal{F}(\widehat{G})$ is known to be a commutative fusion rule algebra with 
unit $Ch(\pi_0) = \pi_0$. 
Let $G_0$ be a closed subgroup of $G$. 
We write $\widehat{G_0} = \{\tau_0, \tau_1, \dots, \tau_n, \dots\}$,  
where $\tau_0$ is the trivial representation of $G_0$.

\medskip

The followings are well-known facts. 

\medskip
\noindent
{\bf Lemma 3.1  } For a compact group $G$, the followings hold.
\begin{itemize}
\item[(1)] $Ch(\pi_i \otimes \pi_j) = Ch(\pi_i) Ch(\pi_j)$ 
for $\pi_i, \pi_j \in  {\rm Rep}^{\rm f}(G)$. 
\item[(2)] $Ch({\rm res}_{G_0}^G \pi) = {\rm res}_{G_0}^G Ch(\pi)$ 
for $\pi \in {\rm Rep}^{\rm f}(G)$. 
\end{itemize}

\medskip
\noindent
{\bf Lemma 3.2 }  For a compact group $G$ and its closed subgroup $G_0$, the 
followings hold. 
\begin{itemize}
\item[(1)] [Character formula] (refer to Hirai [H])\ For $\tau \in {\rm Rep}^{\rm f}(G_0)$, 
$$
Ch({\rm ind}_{G_0}^G \tau)(g) = [G : G_0]\int_G ch(\tau)(sgs^{-1})  
1_{G_0}(sgs^{-1}) \omega_G(ds). 
$$
\item[(2)] [Frobenius' reciprocity theorem] (refer to Folland [F]) \ For $\tau \in \widehat{G_0}$ 
and   for $\pi \in \widehat{G}$,  
$$ 
[{\rm ind}_{G_0}^G \tau : \pi ] = [\tau : {\rm res}_{G_0}^G \pi],  
$$
where $[ \ \ : \ \ ]$ denotes the multiplicity of representations. 
\end{itemize}

\medskip
\noindent
{\bf Lemma 3.3 } For $\tau \in {\rm Rep}^{\rm f}(G_0)$,
$$ 
Ch({\rm ind}_{G_0}^G \tau) = [G : G_0]{\rm ind}_{G_0}^G Ch(\tau).
$$

\medskip
\noindent
{\bf Proof } By the theorem 3.6 in [HKY], 
$$
ch({\rm ind}_{G_0}^G \tau) = {\rm ind}_{G_0}^G ch(\tau)
$$
where 
$$
ch(\tau) := \frac{1}{\dim \tau}Ch(\tau)~~{\rm and}~~
{\rm ind}_{G_0}^G ch(\tau) 
:= \int_{G} ch(\tau)(sgs^{-1}) 1_{G_0}(sgs^{-1})\omega_G(ds). 
$$
Hence 
\begin{align*}
Ch({\rm ind}_{G_0}^G \tau)(g) &= 
[G : G_0]\int_{G}Ch(\tau)(sgs^{-1}) 1_{G_0}(sgs^{-1})\omega_G(ds)\\
&= [G : G_0] {\rm ind}_{G_0}^G Ch(\tau). 
\end{align*}

\medskip
\noindent
{\bf Definition } On the set 
$$
\mathcal{F}(\widehat{G} \cup \widehat{G_0})
:= \{(Ch(\pi), \circ), (Ch(\tau), \bullet) : \pi \in \widehat{G}, \tau \in \widehat{G_0}\}  
$$
we define a convolution $*$ as follows. 
For $\pi_i, \pi_j, \pi \in \widehat{G}$ and 
$\tau_i, \tau_j, \tau \in \widehat{G_0}$, 
\begin{align*}
&(Ch(\pi_i), \circ) * (Ch(\pi_j), \circ) := (Ch(\pi_i)Ch(\pi_j), \circ),\\
&(Ch(\pi), \circ) * (Ch(\tau), \bullet) := (Ch({\rm res}_{G_0}^G \pi) Ch(\tau), \bullet),\\
&(Ch(\tau), \bullet) * (Ch(\pi), \circ) := (Ch(\tau)Ch({\rm res}_{G_0}^G \pi), \bullet),\\
&(Ch(\tau_i), \bullet) * (Ch(\tau_j), \bullet) 
:= (Ch({\rm ind}_{G_0}^G(\tau_i \otimes \tau_j)), \circ).
\end{align*}
We want to check the associativity relations of the convolution in the following cases. 
Whenever reference to a particular representation $\pi$ is not needed, 
we abbreviate  $(Ch(\pi), \circ)$ by $\circ$ and $(Ch(\tau), \bullet)$ 
by $\bullet$. Hence our task will be to verify the subsequent formulae: 
\begin{align*}
&(A1)~(\circ * \circ) * \circ = \circ * (\circ * \circ), \\
&(A2)~(\bullet * \circ) * \circ = \bullet * (\circ * \circ),\\
&(A3)~(\bullet * \bullet) * \circ = \bullet * (\bullet * \circ)~{\rm and}\\
&(A4)~(\bullet * \bullet) * \bullet = \bullet * (\bullet * \bullet). 
\end{align*}

\medskip
\noindent
{\bf Lemma 3.4 } The equalities $(A1)$, $(A2)$ and $(A3)$ hold 
without further assumptions.  
For $\pi_i, \pi_j, \pi_k, \pi \in \widehat{G}$ and 
$\tau_i, \tau_j, \tau \in \widehat{G_0}$, 
\begin{align*}
&(A1)~((Ch(\pi_i), \circ) * (Ch(\pi_j), \circ)) * (Ch(\pi_k), \circ) 
= (Ch(\pi_i), \circ) * ((Ch(\pi_j), \circ) * (Ch(\pi_k), \circ)). \\
&(A2)~ ((Ch(\tau), \bullet) * (Ch(\pi_i), \circ)) * (Ch(\pi_j), \circ)
=  (Ch(\tau), \bullet) * ((Ch(\pi_i), \circ) * (Ch(\pi_j), \circ)).\\
&(A3)~((Ch(\tau_i), \bullet) * (Ch(\tau_j), \bullet)) * (Ch(\pi), \circ)
=(Ch(\tau_i), \bullet) * ((Ch(\tau_j), \bullet) * (Ch(\pi), \circ)).
\end{align*}

\medskip
\noindent
{\bf Proof } $(A1)$ is clear because $\mathcal{F}(\widehat{G})$ is a 
fusion rule algebra.

\medskip

$(A2)$ For $\tau \in \widehat{G_0}$ and $\pi_i, \pi_j \in \widehat{G}$, 
\begin{align*}
&((Ch(\tau), \bullet) * (Ch(\pi_i), \circ)) * (Ch(\pi_j), \circ)\\
=& (Ch(\tau) Ch({\rm res}_{G_0}^G \pi_i), \bullet) * (Ch(\pi_j), \circ)\\
=& (Ch(\tau) Ch({\rm res}_{G_0}^G \pi_i) Ch({\rm res}_{G_0}^G \pi_j) , \bullet)\\
=& (Ch(\tau) ({\rm res}_{G_0}^G Ch(\pi_i)) ({\rm res}_{G_0}^G Ch(\pi_j)) , \bullet). 
\end{align*}
On the other hand, 
\begin{align*}
&(Ch(\tau), \bullet) * ((Ch(\pi_i), \circ) * (Ch(\pi_j), \circ))\\
=& (Ch(\tau), \bullet)) * ((Ch(\pi_i) Ch(\pi_j), \circ)\\
=& (Ch(\tau) {\rm res}_{G_0}^G (Ch(\pi_i) Ch(\pi_j)), \bullet) \\
=& (Ch(\tau) ({\rm res}_{G_0}^G Ch(\pi_i)) ({\rm res}_{G_0}^G Ch(\pi_j)), \bullet).
\end{align*}

\medskip

$(A3)$ For $\tau_i, \tau_j \in \widehat{G_0}$ and $\pi \in \widehat{G}$, 
\begin{align*}
&((Ch(\tau_i), \bullet) * (Ch(\tau_j), \bullet)) * (Ch(\pi), \circ)\\
=& (Ch({\rm ind}_{G_0}^G(\tau_i \otimes \tau_j)), \circ) * (Ch(\pi), \circ)\\
=& (Ch({\rm ind}_{G_0}^G(\tau_i \otimes \tau_j)) Ch(\pi), \circ). 
\end{align*}
For every $g \in G$, 
\begin{align*}
&(Ch({\rm ind}_{G_0}^G(\tau_i \otimes \tau_j)) Ch(\pi))(g) \\
=& Ch({\rm ind}_{G_0}^G(\tau_i \otimes \tau_j)(g) Ch(\pi) (g)\\
=& [G : G_0]\int_G Ch(\tau_i \otimes \tau_j)(sgs^{-1}) 
1_{G_0}(sgs^{-1}) \omega_G(ds)  Ch(\pi) (g) \\
= & [G : G_0]\int_G (Ch(\tau_i) Ch(\tau_j))(sgs^{-1}) Ch(\pi)(g)
1_{G_0}(sgs^{-1}) \omega_G(ds)  \\
= & [G : G_0]\int_G Ch(\tau_i)(sgs^{-1}) Ch(\tau_j)(sgs^{-1})  Ch(\pi) (sgs^{-1})
1_{G_0}(sgs^{-1})\omega_G(ds). 
\end{align*}
On the other hand, 
\begin{align*}
&(Ch(\tau_i), \bullet) * ((Ch(\tau_j), \bullet) * (Ch(\pi), \circ))\\
=& (Ch(\tau_i), \bullet) * (Ch(\tau_j)Ch({\rm res}_{G_0}^G \pi), \bullet)\\
=& (Ch({\rm ind}_{G_0}^G (\tau_i \otimes \tau_j \otimes {\rm res}_{G_0}^G \pi)), \circ).   
\end{align*}
For each $g \in G$, 
\begin{align*}
&Ch({\rm ind}_{G_0}^G (\tau_i \otimes \tau_j \otimes {\rm res}_{G_0}^G \pi))(g) \\
= & [G : G_0]\int_G Ch(\tau_i \otimes \tau_j \otimes {\rm res}_{G_0}^G \pi)(sgs^{-1}) 
1_{G_0}(sgs^{-1}) \omega_G(ds)\\
= & [G : G_0]\int_G (Ch(\tau_i)Ch(\tau_j)Ch({\rm res}_{G_0}^G \pi))(sgs^{-1}) 
1_{G_0}(sgs^{-1}) \omega_G(ds)\\
= & [G : G_0]\int_G Ch(\tau_i)(sgs^{-1})Ch(\tau_j)(sgs^{-1})Ch(\pi)(sgs^{-1}) 
1_{G_0}(sgs^{-1}) \omega_G(ds). 
\end{align*}
\hspace{\fill}$\Box$

\medskip
\noindent
{\bf Definition ([HKTY1])} Let $(G, G_0)$ be a pair of consisting of a compact group $G$ and 
a closed subgroup $G_0$ of $G$. 
For $g \in G_0$
$$
X(g):= \{ s\in G : sgs^{-1} \in G_0 \}.
$$ 
We call $(G, G_0)$ an  {\it  admissible pair}   
if for any $\tau \in \widehat{G_0}$, any $g \in G_0$ and any $s \in X(g)$, 
$$
ch(\tau)(sgs^{-1}) = ch(\tau)(g)
$$
holds.  

\medskip
\noindent
{\bf Remark } It is clear that $Ch(\tau)(sgs^{-1}) = Ch(\tau)(g)$ holds, 
if $(G, G_0)$ is an admissible pair. 

\medskip
\noindent
{\bf Lemma 3.5 } If a compact group $G$ together with a subgroup $G_0$ 
of $G$ with $[G:G_0] < + \infty$ forms an admissible pair, then the associativity relation $(A4)$ holds.

\medskip
\noindent
{\bf Proof } Assume that $(G, G_0)$ is an admissible pair. For $\tau_i, \tau_j, \tau_k \in \widehat{G_0}$ 
and $g \in G_0$ 
\begin{align*}
&((Ch(\tau_i), \bullet) * (Ch(\tau_j), \bullet)) * (Ch(\tau_k), \bullet)\\
=&(Ch({\rm ind}_{G_0}^G (\tau_i \otimes \tau_j)), \circ) * (Ch(\tau_k), \bullet)\\
=&([G : G_0]{\rm ind}_{G_0}^G Ch(\tau_i \otimes \tau_j), \circ) * (Ch(\tau_k), \bullet)\\
=&([G : G_0]({\rm res}_{G_0}^G({\rm ind}_{G_0}^G (Ch(\tau_i) Ch(\tau_j))))Ch(\tau_k), \bullet). 
\end{align*}
For $g \in G_0$, 
\begin{align*}
&({\rm ind}_{G_0}^G (Ch(\tau_i) Ch(\tau_j)))(g) Ch(\tau_k))(g)\\
=& \left([G : G_0]\int_{G} Ch(\tau_i)(sgs^{-1}) Ch(\tau_j)(sgs^{-1}) 1_{G_0}(sgs^{-1})  
\omega_G(ds)\right) Ch(\tau_k)(g)\\
=& \left([G : G_0]\int_{G} Ch(\tau_i)(g) Ch(\tau_j)(g) 1_{G_0}(sgs^{-1})  
\omega_G(ds)\right) Ch(\tau_k)(g)\\
=& \left([G : G_0]\int_{G} 1_{G_0}(sgs^{-1})  
\omega_G(ds)\right) Ch(\tau_i)(g) Ch(\tau_j)(g) Ch(\tau_k)(g). 
\end{align*}
This implies the associativity relation $(A4)$. 
\hspace{\fill}$\Box$

\medskip
\noindent
{\bf Lemma 3.6 } If  the associativity relation $(A4)$ holds  
for a compact group $G$ and a subgroup $G_0$ of $G$ with $[G:G_0] < + \infty$, 
then $(G,G_0)$  is an admissible pair. 

\medskip
\noindent
{\bf Proof }
Assume that the associativity relation $(A4)$ holds. Let $\tau_0$ be the trivial 
representation of $\widehat{G_0}$. For $\tau \in \widehat{G_0}$ the associativity 
relation  
$$
((Ch(\tau_0), \bullet) * (Ch(\tau_0), \bullet)) * (Ch(\tau), \bullet) 
= (Ch(\tau_0), \bullet) * ((Ch(\tau_0), \bullet) * (Ch(\tau), \bullet)) 
$$
holds. 
\begin{align*}
((Ch(\tau_0), \bullet) * (Ch(\tau_0), \bullet)) * (Ch(\tau), \bullet)
= (Ch({\rm res}_{G_0}^G ({\rm ind}_{G_0}^G  \tau_0 )) Ch(\tau),\bullet) 
\end{align*}
and
$$
(Ch(\tau_0), \bullet) * ((Ch(\tau_0), \bullet) * (Ch(\tau), \bullet))
= (Ch({\rm res}_{G_0}^G ({\rm ind}_{G_0}^G \tau)), \bullet). 
$$
Then for $g \in G_0$
$$
Ch({\rm ind}_{G_0}^G \tau_0)(g) Ch(\tau)(g) = Ch({\rm ind}_{G_0}^G \tau)(g). 
$$
Now  
$$
Ch({\rm ind}_{G_0}^{G} \tau_0)(g) \geq 1. 
$$
Indeed by the character formula 
\begin{align*}
Ch({\rm ind}_{G_0}^{G} \tau_0)(g)  
&= [G : G_0]\int_{G} Ch(\tau_0)(sgs^{-1}) 1_{G_0} (sgs^{-1})  \omega_{G}(ds)\\ 
&= [G : G_0] \int_{G} 1_{G_0}(sgs^{-1}) \omega_G(ds)\\
&=[G : G_0]\omega_{G}(X(g)). 
\end{align*}
Since $X(g) \supset G_0$, we see that $\omega_{G}(X(g)) \geq \omega_{G}(G_0)$. 
By the assumption  $[G :G_0] < + \infty$, we obtain $\omega_{G}(G_0) = 1/[G :G_0] > 0$. 
Therefore 
$$
[G : G_0]\omega_{G}(X(g)) \geq [G : G_0]\omega_{G}(G_0) = 1. 
$$
Hence 
$$
Ch(\tau)(g) = (Ch({\rm ind}_{G_0}^G \tau_0)(g))^{-1} Ch({\rm ind}_{G_0}^G \tau)(g). 
$$
For $s \in X(g)$
\begin{align*}
Ch(\tau)(sgs^{-1}) 
&= (Ch({\rm ind}_{G_0}^{G} \tau_0)(sgs^{-1}))^{-1}  Ch({\rm ind}_{G_0}^{G} \tau)(sgs^{-1})\\
&= (Ch({\rm ind}_{G_0}^{G} \tau_0)(g))^{-1}  Ch({\rm ind}_{G_0}^{G} \tau)(g)\\
&= Ch(\tau)(g).  
\end{align*}
Then $(G, G_0)$ is an admissible pair. 
\hspace{\fill}$\Box$

\bigskip
\noindent
{\bf Theorem } Let $G_0$ be a closed subgroup of a compact group $G$ such that 
$[G:G_0] < + \infty$. Then $\mathcal{F}(\widehat{G} \cup \widehat{G_0})$ is a 
fusion rule algebra if and only if $(G, G_0)$ is an admissible pair.

\medskip
\noindent
{\bf Proof } The associativity relations $(A1)$, $(A2)$ and $(A3)$ are 
a consequence of Lemma 3.4, and $(A4)$ holds if and only if 
$(G, G_0)$ is an admissible pair by Lemma 3.5 and  Lemma 3.6. 
It is easy to check the remaining axioms of a fusion rule algebra  
for $\mathcal{F}(\widehat{G} \cup \widehat{G_0})$. 
The desired conclusion follows. 
\hspace{\fill}$\Box$

\medskip
\noindent
{\bf Remark} 
\begin{itemize}
\item[(1)] The above $\mathcal{F}(\widehat{G} \cup \widehat{G_0})$ is a fusion rule algebra 
such that the sequence: 
$$
1 \longrightarrow \mathcal{F}(\widehat{G}) \longrightarrow \mathcal{F}(\widehat{G} \cup \widehat{G_0}) 
\longrightarrow \mathbb{Z}_2 \longrightarrow 1
$$
is exact. 
\item[(2)] If $G_0 = G$, then $\mathcal{F}(\widehat{G} \cup \widehat{G_0})$ is the 
fusion rule algebra $\mathcal{F}(\widehat{G}) \times \mathbb{Z}_2$. 
\item[(3)] If $G$ is a finite group and $G_0 = \{e\}$, then  
$\mathcal{F}(\widehat{G} \cup \widehat{G_0})$ is the fusion rule algebra join 
$\mathcal{F}(\widehat{G}) \vee_F \mathbb{Z}_2$. 
\item[(4)] $\mathcal{K}_d(\mathcal{F}(\widehat{G} \cup \widehat{G_0})) 
\cong \mathcal{K}(\widehat{G} \cup \widehat{G_0})$, where  $\mathcal{K}(\widehat{G} \cup \widehat{G_0})$ 
is a hypergroup in the sense that of [HKTY1]. 
\end{itemize}

\bigskip

The following lemmas (Lemma 3.7, 3.9, 3.11, 3.13, 3.15 and 3.16) are quoted from [HKTY1].  
The following corollaries (Corollary 3.8, 3.10, 3.12 and 3.14) are rewrote version of 
``hypergroup'' in [HKTY1] to ``fusion rule algebra''. 
The proofs of them are almost same as in the case of [HKTY1], hence 
we omit the proofs. Regarding proofs of them, refer to [HKTY1].

\medskip
\noindent
{\bf Lemma 3.7 ([HKTY1])} If $G$ is a compact Abelian group and $G_0$ a closed subgroup of $G$ 
with $[G: G_0] < + \infty$. Then $(G, G_0)$ is always an admissible pair. 

%\medskip
%\noindent
%{\bf Proof } The desired assertion clearly follows from the fact that $sgs^{-1} = g$ 
%for $g \in G_0$ and $s \in G$. 
%\hspace{\fill}$\Box$

\bigskip
\noindent
{\bf Corollary 3.8 } Let $G$ be a compact Abelian group and $G_0$ a closed subgroup of $G$ 
with $[G: G_0] < + \infty$. Then $\mathcal{F}(\widehat{G} \cup \widehat{G_0})$ 
is a fusion rule algebra. 

\if0
\medskip
\noindent
{\bf Proof } This assertion follows directly from the Theorem  and Lemma 3.7. 
\hspace{\fill}$\Box$
\fi

\medskip
\noindent
{\bf Lemma 3.9 ([HKTY1])} If for each $\tau \in \widehat{G_0}$ there exists a representation 
$\tilde{\tau}$ of $G$ such that ${\rm res}_{G_0}^G \tilde{\tau} = \tau$, then 
$(G, G_0)$ is an admissible pair.

%\medskip
%\noindent
%{\bf Proof }  For $\tau \in \widehat{G_0}$, $g \in G_0$ and $s \in X(g)$, 
%$$
%ch(\tau)(sgs^{-1}) = ch(\tilde{\tau})(sgs^{-1}) = ch(\tilde{\tau})(g) = ch(\tau)(g). 
%$$
%\hspace{\fill}$\Box$

\bigskip
\noindent
{\bf Corollary 3.10 } Let $G$ be a semi-direct product group 
$H \rtimes_\alpha G_0$,  where $H$ is a finite group and $G_0$ is a 
finite group. Then 
$\mathcal{F}(\widehat{G} \cup \widehat{G_0})$ is a fusion rule algebra.

\if0
\medskip
\noindent
{\bf Proof } For $\tau \in \widehat{G_0}$, put 
$$
\tilde{\tau}((h,g)) = \tau(g)
$$
for $(h,g) \in H \rtimes_\alpha G_0 = G$. Then $\tilde{\tau}$ is a finite 
dimensional representation of $G$ and ${\rm res}_{G_0}^G \tilde{\tau} = \tau.$ 
By the Theorem  and Lemma 3.9 we arrive at the desired conclusion. 

\hspace{\fill}$\Box$
\fi

\medskip
\noindent
{\bf Lemma 3.11 ([HKTY1])} If for $g \in G_0$ and $s \in X(g)$ there exists $t \in G_0$ 
such that $tgt^{-1} = sgs^{-1}$, then $(G, G_0)$ is an admissible pair. 

%\medskip
%\noindent
%{\bf Proof } For $ \tau \in \widehat{G_0}$, $g \in G_0$ and $s \in X(g)$, 
%$$
%ch(\tau)(sgs^{-1}) = ch(\tau)(tgt^{-1}) = ch(\tau)(g). 
%$$  
%\hspace{\fill}$\Box$

\bigskip

Let $S_n$ be the symmetric group of degree $n$.

\medskip
\noindent
{\bf Corollary 3.12 } $\mathcal{F}(\widehat{S_n} \cup \widehat{S_{n-1}})$ $(n \geq 2)$ is a 
fusion rule algebra.

\if0 
\medskip
\noindent
{\bf Proof }  For  $g \in S_{n-1}$, $s \in X(g)$ such that $s^{-1}(n) = a$, 
$$
sgs^{-1}(n) = s(g(s^{-1}(n))) = s(g(a)). 
$$
Since $sgs^{-1} \in S_{n-1}$, $sgs^{-1}(n) = n$ and $s(g(a)) = n$. 
Then $g(a) = s^{-1}(n) = a$. 

Put $t = s s_1$ where $s_1$ is a transposed permutation $(a, n)$. Then 
we  see that 
$$
tgt^{-1} = sgs^{-1} 
$$
by the fact that for $b$ such that $b \not = a$, $s^{-1}(b) \not = a$ and 
$g(s^{-1}(b)) \not = a$ hold. By the Theorem  and Lemma 3.11 we get the desired conclusion. 
\hspace{\fill}$\Box$
\fi

\medskip
\noindent
{\bf Lemma 3.13 ([HKTY1])} Let $G_0$ and $G_1$ be  closed subgroups of $G$ such that 
$G_0 \subset G_1 \subset G$. If $(G_1, G_0)$ and $(G, G_1)$ are admissible pairs, 
then $(G, G_0)$ is an admissible pair.

\if0
\medskip
\noindent
{\bf Proof } Since $(G_1, G_0)$ is an admissible pair, for $\tau \in \widehat{G_0}$ and 
$g \in G_0$, 
\begin{align*}
ch({\rm ind}_{G_0}^{G_1} \tau)(g)  
%
%
&= \int_{G_1} ch(\tau)(sgs^{-1}) 1_{G_0} (sgs^{-1}) d \omega_{G_1}(s)\\
%
%
&= \int_{G_1} ch(\tau)(g) 1_{G_0} (sgs^{-1}) d \omega_{G_1}(s)\\
%
%
&= \int_{G_1}  1_{G_0} (sgs^{-1}) d \omega_{G_1}(s) ch(\tau)(g)\\
%
%
&= ch({\rm ind}_{G_0}^{G_1} \tau_0)(g) ch(\tau)(g), 
\end{align*}
where $\tau_0$ is the trivial representation of $G_0$. 
Since 
$$
ch({\rm res}_{G_0}^{G_1} \tau_0)(g) \geq \omega_{G_1}(G_0) >0,
$$  
we see that
$$
ch(\tau)(g) =  (ch({\rm ind}_{G_0}^{G_1} \tau_0)(g))^{-1}  ch({\rm ind}_{G_0}^{G_1} \tau)(g). 
$$
Since  $(G, G_1)$ is an admissible pair and ${\rm ind}_{G_0}^{G_1} \tau_0 \in {\rm Rep^f}(G_1)$,   
${\rm ind}_{G_0}^{G_1} \tau \in {\rm Rep^f}(G_1)$, we have for $g \in G_0 \subset G_1$ and 
$s \in X(g)$, 
$$
ch({\rm ind}_{G_0}^{G_1} \tau_0)(sgs^{-1}) = ch({\rm ind}_{G_0}^{G_1} \tau_0)(g)
$$
and 
$$
ch({\rm ind}_{G_0}^{G_1} \tau)(sgs^{-1}) = ch({\rm ind}_{G_0}^{G_1} \tau)(g). 
$$
Then we obtain 
\begin{align*}
ch(\tau)(sgs^{-1}) 
%
%
&= (ch({\rm ind}_{G_0}^{G_1} \tau_0)(sgs^{-1}))^{-1}  ch({\rm ind}_{G_0}^{G} \tau)(sgs^{-1})\\
%
%
&= (ch({\rm ind}_{G_0}^{G_1} \tau_0)(g))^{-1}  ch({\rm ind}_{G_0}^{G_1} \tau)(g)\\
%
%
&= ch(\tau)(g) 
\end{align*}
\hspace{\fill}$\Box$
\fi

\bigskip
\noindent
{\bf Corollary 3.14 } For natural numbers $m$ and $n$ such that 
$m > n \geq 1$,  $\mathcal{F}(\widehat{S_m} \cup \widehat{S_n})$ is a 
fusion rule algebra.

\if0
\medskip
\noindent
{\bf Proof } This statement follows from the Theorem  and Lemma 3.13 combined with 
Corollary 3.12.  
\hspace{\fill}$\Box$
\fi

\bigskip
Let $G_0$ be a closed normal subgroup of $G$. Then the coadjoint action 
$\widehat{\alpha}$ of $G$ on $\widehat{G_0}$ is defined by 
$$
\widehat{\alpha}_s (\tau)(g) := \tau(sgs^{-1})
$$
for $\tau \in \widehat{G_0}$, $g \in G_0$ and $s \in G$. 
If $\widehat{\alpha}_s = id $ for all $s \in G$, we say that  $\widehat{\alpha}$ is trivial.

\medskip
\noindent
{\bf Lemma 3.15 ([HKTY1])} Let $G_0$ be a closed normal subgroup of $G$. The 
pair $(G, G_0)$ is an admissible pair if and only if the coadjoint action 
$\widehat{\alpha}$ is trivial. 

\if0
\medskip
\noindent
{\bf Proof  } Assume that $(G, G_0)$ is an admissible pair. For $g \in G_0$ it 
is clear that $X(g) = G$. Then for $\tau \in \widehat{G_0}$
$$
ch(\tau)(sgs^{-1}) = ch(\tau)(g) 
$$ 
for all $s \in G$. This implies that 
$$
ch(\widehat{\alpha}_s(\tau))(g) = ch(\tau)(g) 
$$
for all $g \in G_0$. Hence we obtain 
$$
\hat{\alpha}_s(\tau) \cong \tau
$$
for $\tau \in \widehat{G_0}$. In fact $\hat{\alpha}_s$ is the identity on $\widehat{G_0}$ which 
means that $\hat{\alpha}$ is trivial. 

The converse is clear.  
\hspace{\fill}$\Box$
\fi

\bigskip
\noindent
{\bf Lemma 3.16 ([HKTY1])} Let $G_0$ be a closed normal commutative subgroup of $G$. 
The pair $(G, G_0)$ is an admissible pair if and only if $G \cong G_0 \times (G/{G_0})$.

\if0
\bigskip
\noindent
{\bf Proof  } Assume that $(G, G_0)$ is an admissible pair and $sgs^{-1} \not = g$ for 
$g \in G_0$ and $s \in X(g) = G$. Since $\widehat{G_0}$ separates $G_0$, there exists 
$\tau \in \widehat{G_0}$ such that
$$
\tau(sgs^{-1}) \not = \tau(g). 
$$ 
This contradicts the assumption that $(G, G_0)$ is an admissible pair. 
But then for $g \in G_0$ and $s \in G$
$$
sgs^{-1} = g, 
$$
namely
$$
sg = gs
$$
holds. For $g_1, g_2 \in G_0$ and $s_1, s_2 \in G$
$$
(g_1 s_1) (g_2 s_2) = g_1 (s_1 g_2 ) s_2 = g_1(g_2 s_1) s_2 = (g_1 g_2)(s_1 s_2).
$$
This implies that $G \cong G_0 \times (G/G_0)$. 

The converse is clear by Lemma 3.15. 
\hspace{\fill}$\Box$
\fi

\bigskip
\noindent
{\bf Corollary 3.17 } Let $G$ be a semi-direct product group $H \rtimes_\alpha G_0$ 
where  $H$ is a compact Abelian group and $G_0$ is a finite group. 
$\mathcal{F}(\widehat{G} \cup \widehat{H})$ is a fusion rule algebra if and only if 
the action $\alpha$ is trivial, i.e. $G = H \times G_0$.

\if0
\medskip
\noindent
{\bf Proof } We note that $H$ is a closed normal subgroup of $G$ and $G/H \cong G_0$. 
Then the assertion follows from the Theorem together with Lemma 3.16. 
\hspace{\fill}$\Box$
\fi

\section{Examples}

\medskip

Associated with a pair $(G, G_0)$ of finite groups such that $G \supset G_0$, 
we obtain a certain finite graph $D(\widehat{G} \cup \widehat{G_0})$ by Frobenius' 
reciprocity theorem. The set of vertices is  
$$
 \{(\pi, \circ), (\tau, \bullet) : \pi \in \widehat{G}, \tau \in \widehat{G_0}\}  
$$
and the edge between $(\pi, \circ)$ and $(\tau, \bullet)$ is given by the 
multiplicity 
$$
m_{\pi,\tau} := [{\rm ind}_{G_0}^G(\tau) : \pi ] = [\tau : {\rm res}_{G_0}^G \pi] \neq 0.
$$
We call this graph $D(\widehat{G} \cup \widehat{G_0})$ a Frobenius diagram. 
Frobenius diagrams $D(\widehat{G} \cup \widehat{G_0})$ sometimes appear 
as Dynkin diagrams and sometimes as Coxter graphs ([GHJ]). V. S. Sunder and N. J. Wildberger 
constructed in [SW] fusion rule algebras $\mathcal{F}(D)$ and hypergroups 
$\mathcal{K}(D)$ associated with certain Dynkin diagrams of type $A_n$,  $D_{2n}$ and so on. 
We give some examples of fusion rule algebras $\mathcal{F}(\widehat{G} \cup \widehat{G_0})$ 
which are compatible with Frobenius diagrams $D(\widehat{G} \cup \widehat{G_0})$.

\bigskip
\noindent
{\bf 4.1 } \ The case that 
$G = \mathbb{Z}_2 = \{e,g\}$ $(g^2 = e)$ and $G_0 = \{e\}$. 
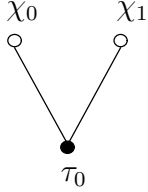
\begin{center}
%WinTpicVersion4.30a
{\unitlength 0.1in%
\begin{picture}(  6.2400,  8.5900)( 19.8000,-17.0000)%
% LINE 2 0 3 0 Black White
% 4 2021 1109 2295 1611 2302 1604 2570 1109
% 
\special{pn 8}%
\special{pa 2021 1109}%
\special{pa 2295 1611}%
\special{fp}%
\special{pa 2302 1604}%
\special{pa 2570 1109}%
\special{fp}%
% STR 2 0 3 0 Black White
% 4 1980 909 1980 971 2 0 0 0
% $\chi_0$
\put(19.8000,-9.7100){\makebox(0,0)[lb]{$\chi_0$}}%
% STR 2 0 3 0 Black White
% 4 2549 909 2549 971 2 0 0 0
% $\chi_1$
\put(25.4900,-9.7100){\makebox(0,0)[lb]{$\chi_1$}}%
% STR 2 0 3 0 Black White
% 4 2261 1767 2261 1830 2 0 0 0
% $\tau_0$
\put(22.6100,-18.3000){\makebox(0,0)[lb]{$\tau_0$}}%
% ELLIPSE 2 0 3 0 Black White
% 4 2570 1072 2604 1104 2572 1040 2572 1040
% 
\special{pn 8}%
\special{ar 2570 1072 34 32  0.0000000  6.2831853}%
% ELLIPSE 2 0 3 0 Black White
% 4 2021 1072 2056 1104 2023 1040 2023 1040
% 
\special{pn 8}%
\special{ar 2021 1072 35 32  0.0000000  6.2831853}%
% ELLIPSE 2 0 0 0 Black Black
% 4 2295 1629 2330 1661 2297 1598 2297 1598
% 
\special{sh 1.000}%
\special{ia 2295 1629 35 32  0.0000000  6.2831853}%
\special{pn 8}%
\special{ar 2295 1629 35 32  0.0000000  6.2831853}%
\end{picture}}%
\end{center}
\bigskip
$\mathcal{F}(\widehat{G} \cup \widehat{G_0}) = \{(Ch(\chi_0), \circ), (Ch(\chi_1), \circ), 
(Ch(\tau_0), \bullet)\}$. Put $\gamma_0 = (Ch(\chi_0), \circ)$, 
$\gamma_1 = (Ch(\chi_1), \circ)$ and $\rho_0 = (Ch(\tau_0), \bullet)$. Then 
the structure equations are
\begin{align*}
\gamma_1 \gamma_1 = \gamma_0,~~~\rho_0 \rho_0 = \gamma_0 + \gamma_1,~~~
\gamma_1 \rho_0 = \rho_0. 
\end{align*} 
Then, the values of the dimension function of 
$\mathcal{F}(\widehat{G} \cup \widehat{G_0})$ are 
\begin{align*}
d(\gamma_0) = 1,~~~d(\gamma_1) = 1,~~~d(\rho_0) = \sqrt{2}. 
\end{align*}
Put 
\begin{align*}
\tilde{\gamma_0} := \frac{1}{d(\gamma_0)}\gamma_0,~~~
\tilde{\gamma_1} := \frac{1}{d(\gamma_1)}\gamma_1,~~~
\tilde{\rho_0} := \frac{1}{d(\rho_0)}\rho_0.
\end{align*}
The set $\{\tilde{\gamma_0},\tilde{\gamma_1}, \tilde{\rho_0}\}$ becomes a hypergroup, and 
isomorphic to the hypergroup $\mathcal{K}(\widehat{G} \cup \widehat{G_0})$ of 
[HKTY1, example 4.1]

\bigskip
\noindent
{\bf 4.2 } \ The case that 
$G = \mathbb{Z}_3 = \{e,g, g^2\}$ $(g^3 = e)$ and $G_0 = \{e\}$. 
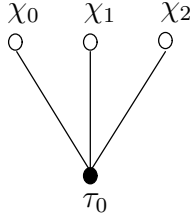
\begin{center}
%WinTpicVersion4.30a
{\unitlength 0.1in%
\begin{picture}(  8.5300,  9.9300)( 19.8000,-18.9000)%
% STR 2 0 3 0 Black White
% 4 1980 953 1980 1035 2 0 0 0
% $\chi_0$
\put(19.8000,-10.3500){\makebox(0,0)[lb]{$\chi_0$}}%
% STR 2 0 3 0 Black White
% 4 2375 953 2375 1035 2 0 0 0
% $\chi_1$
\put(23.7500,-10.3500){\makebox(0,0)[lb]{$\chi_1$}}%
% STR 2 0 3 0 Black White
% 4 2368 1938 2368 2020 2 0 0 0
% $\tau_0$
\put(23.6800,-20.2000){\makebox(0,0)[lb]{$\tau_0$}}%
% ELLIPSE 2 0 3 0 Black White
% 4 2409 1133 2444 1175 2411 1092 2411 1092
% 
\special{pn 8}%
\special{ar 2409 1133 35 42  0.0000000  6.2831853}%
% ELLIPSE 2 0 3 0 Black White
% 4 2021 1133 2056 1175 2023 1092 2023 1092
% 
\special{pn 8}%
\special{ar 2021 1133 35 42  0.0000000  6.2831853}%
% ELLIPSE 2 0 0 0 Black Black
% 4 2409 1831 2444 1873 2411 1789 2411 1789
% 
\special{sh 1.000}%
\special{ia 2409 1831 35 42  0.0000000  6.2831853}%
\special{pn 8}%
\special{ar 2409 1831 35 42  0.0000000  6.2831853}%
% ELLIPSE 2 0 3 0 Black White
% 4 2798 1125 2833 1167 2800 1083 2800 1083
% 
\special{pn 8}%
\special{ar 2798 1125 35 42  0.0000000  6.2831853}%
% LINE 2 0 3 0 Black White
% 6 2409 1183 2409 1807 2798 1183 2409 1807 2028 1183 2409 1807
% 
\special{pn 8}%
\special{pa 2409 1183}%
\special{pa 2409 1807}%
\special{fp}%
\special{pa 2798 1183}%
\special{pa 2409 1807}%
\special{fp}%
\special{pa 2028 1183}%
\special{pa 2409 1807}%
\special{fp}%
% STR 2 0 3 0 Black White
% 4 2771 945 2771 1027 2 0 0 0
% $\chi_2$
\put(27.7100,-10.2700){\makebox(0,0)[lb]{$\chi_2$}}%
\end{picture}}%
\end{center}
\bigskip
$\mathcal{F}(\widehat{G} \cup \widehat{G_0}) = \{(Ch(\chi_0), \circ), (Ch(\chi_1), \circ), 
(Ch(\chi_2), \circ), (Ch(\tau_0), \bullet)\}$. Put $\gamma_0 = (Ch(\chi_0), \circ)$, 
$\gamma_1 = (Ch(\chi_1), \circ)$, $\gamma_2 = (Ch(\chi_2), \circ)$ 
and $\rho_0 = (Ch(\tau_0), \bullet)$. Then 
the structure equations are
\begin{align*}
&\gamma_1 \gamma_1 = \gamma_2,~~~\gamma_2 \gamma_2 = \gamma_1,~~~
\gamma_1 \gamma_2 = \gamma_0,\\
&\rho_0 \rho_0 = \gamma_0 + \gamma_1 + \gamma_2,~~~
\gamma_1 \rho_0 = \rho_0,~~~\gamma_2 \rho_0 = \rho_0.  
\end{align*} 
As in the case of example 4.1, the values of the dimension function of 
$\mathcal{F}(\widehat{G} \cup \widehat{G_0})$ are 
\begin{align*}
d(\gamma_0) = 1,~~~d(\gamma_1) = 1,~~~d(\gamma_2)=1,~~~d(\rho_0) = \sqrt{3}. 
\end{align*}
Put 
\begin{align*}
\tilde{\gamma_i} := \frac{1}{d(\gamma_i)}\gamma_i~~~(i=0,1,2),~~~
\tilde{\rho_0} := \frac{1}{d(\rho_0)}\rho_0.
\end{align*}
The set $\{\tilde{\gamma_0},\tilde{\gamma_1}, \tilde{\gamma_2}, \tilde{\rho_0}\}$ becomes a hypergroup, 
and isomorphic to the hypergroup $\mathcal{K}(\widehat{G} \cup \widehat{G_0})$ of 
[HKTY1, example 4.2].

\bigskip
\noindent
{\bf 4.3 } \ The case that $G$ is the symmetric group 
$S_3 = \mathbb{Z}_3 \rtimes_\alpha \mathbb{Z}_2$ of degree 3 and  $G_0 = \mathbb{Z}_2$.
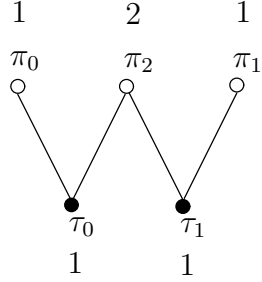
\begin{center}
%WinTpicVersion4.32a
{\unitlength 0.1in%
\begin{picture}(12.1200,13.2000)(19.8000,-19.6000)%
% LINE 2 0 3 0 Black White  
% 4 2022 1146 2300 1702 2307 1695 2578 1146
% 
\special{pn 8}%
\special{pa 2022 1146}%
\special{pa 2300 1702}%
\special{fp}%
\special{pa 2307 1695}%
\special{pa 2578 1146}%
\special{fp}%
% STR 2 0 3 0 Black White  
% 4 1980 951 1980 1020 2 0 0 0
% $\pi_0$
\put(19.8000,-10.2000){\makebox(0,0)[lb]{$\pi_0$}}%
% STR 2 0 3 0 Black White  
% 4 3135 958 3135 1027 2 0 0 0
% $\pi_1$
\put(31.3500,-10.2700){\makebox(0,0)[lb]{$\pi_1$}}%
% STR 2 0 3 0 Black White  
% 4 2286 1813 2286 1883 2 0 0 0
% $\tau_0$
\put(22.8600,-18.8300){\makebox(0,0)[lb]{$\tau_0$}}%
% CIRCLE 2 0 3 0 Black White  
% 4 2585 1104 2587 1068 2587 1068 2587 1068
% 
\special{pn 8}%
\special{ar 2585 1104 36 36 0.0000000 6.2831853}%
% CIRCLE 2 0 3 0 Black White  
% 4 2022 1104 2024 1068 2024 1068 2024 1068
% 
\special{pn 8}%
\special{ar 2022 1104 36 36 0.0000000 6.2831853}%
% CIRCLE 2 0 0 0 Black Black  
% 4 2300 1723 2302 1688 2302 1688 2302 1688
% 
\special{sh 1.000}%
\special{ia 2300 1723 35 35 0.0000000 6.2831853}%
\special{pn 8}%
\special{ar 2300 1723 35 35 0.0000000 6.2831853}%
% LINE 2 0 3 0 Black White  
% 4 2599 1139 2877 1696 2884 1689 3156 1139
% 
\special{pn 8}%
\special{pa 2599 1139}%
\special{pa 2877 1696}%
\special{fp}%
\special{pa 2884 1689}%
\special{pa 3156 1139}%
\special{fp}%
% CIRCLE 2 0 3 0 Black White  
% 4 3156 1097 3158 1061 3158 1061 3158 1061
% 
\special{pn 8}%
\special{ar 3156 1097 36 36 0.0000000 6.2831853}%
% CIRCLE 2 0 0 0 Black Black  
% 4 2877 1731 2879 1695 2879 1695 2879 1695
% 
\special{sh 1.000}%
\special{ia 2877 1731 36 36 0.0000000 6.2831853}%
\special{pn 8}%
\special{ar 2877 1731 36 36 0.0000000 6.2831853}%
% STR 2 0 3 0 Black White  
% 4 2863 1820 2863 1890 2 0 0 0
% $\tau_1$
\put(28.6300,-18.9000){\makebox(0,0)[lb]{$\tau_1$}}%
% STR 2 0 3 0 Black White  
% 4 2571 958 2571 1027 2 0 0 0
% $\pi_2$
\put(25.7100,-10.2700){\makebox(0,0)[lb]{$\pi_2$}}%
% STR 2 0 3 0 Black White  
% 4 1990 670 1990 770 2 0 0 0
% 1
\put(19.9000,-7.7000){\makebox(0,0)[lb]{1}}%
% STR 2 0 3 0 Black White  
% 4 3150 670 3150 770 2 0 0 0
% 1
\put(31.5000,-7.7000){\makebox(0,0)[lb]{1}}%
% STR 2 0 3 0 Black White  
% 4 2580 680 2580 780 2 0 0 0
% 2
\put(25.8000,-7.8000){\makebox(0,0)[lb]{2}}%
% STR 2 0 3 0 Black White  
% 4 2280 1980 2280 2080 2 0 0 0
% 1
\put(22.8000,-20.8000){\makebox(0,0)[lb]{1}}%
% STR 2 0 3 0 Black White  
% 4 2860 1990 2860 2090 2 0 0 0
% 1
\put(28.6000,-20.9000){\makebox(0,0)[lb]{1}}%
\end{picture}}%
\end{center}
\bigskip
$\mathcal{F}(\widehat{G} \cup \widehat{G_0}) 
= \{(Ch(\pi_i), \circ), (Ch(\tau_j), \bullet): \pi_i \in \widehat{G}, \tau_j \in \widehat{G_0}\}$. 
Put $\gamma_i = (Ch(\pi_i), \circ)$ and $\rho_j = (Ch(\tau_j), \bullet)$.  
Then the structure equations are
\begin{align*}
&\gamma_1 \gamma_1 = \gamma_0,~~~
\gamma_2 \gamma_2 = \gamma_0 + \gamma_1 + \gamma_2,~~~
\gamma_1 \gamma_2 = \gamma_2,~~~
\rho_0 \rho_0 = \rho_1 \rho_1 = \gamma_0 + \gamma_2,\\
&\rho_0 \rho_1 = \rho_1 \rho_0 = \gamma_1 + \gamma_2,~~~
\gamma_0 \rho_0 = \rho_0,~~~\gamma_1 \rho_0 = \rho_1,~~~
\gamma_2 \rho_0 = \rho_0 + \rho_1,\\
&\gamma_0 \rho_1 = \rho_1,~~~\gamma_1 \rho_1 = \rho_0,~~~
\gamma_2 \rho_1 = \rho_0 + \rho_1. 
\end{align*}

\medskip
\noindent
{\bf Remark } $\mathcal{F}(\widehat{\mathbb{Z}_2} \cup \widehat{\{e\}}) = 
\mathcal{F}(A_3)$, 
$\mathcal{F}(\widehat{\mathbb{Z}_3} \cup \widehat{\{e\}}) = 
\mathcal{F}(D_4)$ and 
$\mathcal{F}(\widehat{S_3} \cup \widehat{\mathbb{Z}_2}) = \mathcal{F}(A_5)$ 
where $\mathcal{F}(A_3)$, $\mathcal{F}(D_4)$ and 
$\mathcal{F}(A_5)$ are Sunder-Wildberger's fusion rule algebra (\cite{SW}) associated with  
Dynkin diagrams of type $A_3$, $D_4$ and $A_5$ respectively.

\medskip
As in the case of example 4.1, 
$\mathcal{F}(\widehat{G} \cup \widehat{G_0})$ are 
\begin{align*}
d(\gamma_0) = 1,~~~d(\gamma_1) = 1,~~~d(\gamma_2)=2,~~~
d(\rho_0) = \sqrt{3},~~~d(\rho_1) = \sqrt{3}. 
\end{align*}
Put 
\begin{align*}
\tilde{\gamma_i} := \frac{1}{d(\gamma_i)}\gamma_i~~~(i=0,1,2),~~~
\tilde{\rho_j} := \frac{1}{d(\rho_j)}\rho_j~~~(j=0,1).
\end{align*}
The set 
$\{\tilde{\gamma_0},\tilde{\gamma_1}, \tilde{\gamma_2}, \tilde{\rho_0}, \tilde{\rho_1}\}$ 
becomes a hypergroup, and isomorphic to the hypergroup  
of [HKTY1, example 4.3].

\bigskip
\noindent
{\bf 4.4 } \ The case that 
$G = \mathbb{Z}_4 = \{e,g, g^2, g^3\}$ $(g^4 = e)$ and $G_0 = \mathbb{Z}_2 $.
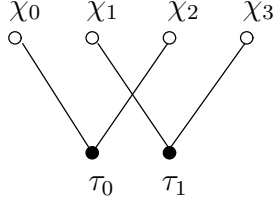
\begin{center}
%WinTpicVersion4.32a
{\unitlength 0.1in%
\begin{picture}(12.6400,9.1000)(5.6800,-18.9000)%
% STR 2 0 3 0 Black White  
% 4 570 1043 570 1110 2 0 0 0
% $\chi_0$
\put(5.7000,-11.1000){\makebox(0,0)[lb]{$\chi_0$}}%
% STR 2 0 3 0 Black White  
% 4 970 1044 970 1110 2 0 0 0
% $\chi_1$
\put(9.7000,-11.1000){\makebox(0,0)[lb]{$\chi_1$}}%
% STR 2 0 3 0 Black White  
% 4 980 1954 980 2020 2 0 0 0
% $\tau_0$
\put(9.8000,-20.2000){\makebox(0,0)[lb]{$\tau_0$}}%
% STR 2 0 3 0 Black White  
% 4 1360 1953 1360 2020 2 0 0 0
% $\tau_1$
\put(13.6000,-20.2000){\makebox(0,0)[lb]{$\tau_1$}}%
% STR 2 0 3 0 Black White  
% 4 1390 1044 1390 1110 2 0 0 0
% $\chi_2$
\put(13.9000,-11.1000){\makebox(0,0)[lb]{$\chi_2$}}%
% CIRCLE 2 0 3 0 Black White  
% 4 1000 1200 1030 1210 1030 1210 1030 1210
% 
\special{pn 8}%
\special{ar 1000 1200 32 32 0.0000000 6.2831853}%
% CIRCLE 2 0 0 0 Black Black  
% 4 1400 1800 1430 1810 1430 1810 1430 1810
% 
\special{sh 1.000}%
\special{ia 1400 1800 32 32 0.0000000 6.2831853}%
\special{pn 8}%
\special{ar 1400 1800 32 32 0.0000000 6.2831853}%
% CIRCLE 2 0 3 0 Black White  
% 4 1400 1200 1430 1210 1430 1210 1430 1210
% 
\special{pn 8}%
\special{ar 1400 1200 32 32 0.0000000 6.2831853}%
% CIRCLE 2 0 0 0 Black Black  
% 4 1000 1800 1030 1810 1030 1810 1030 1810
% 
\special{sh 1.000}%
\special{ia 1000 1800 32 32 0.0000000 6.2831853}%
\special{pn 8}%
\special{ar 1000 1800 32 32 0.0000000 6.2831853}%
% CIRCLE 2 0 3 0 Black Black  
% 4 1800 1200 1830 1210 1830 1210 1830 1210
% 
\special{pn 8}%
\special{ar 1800 1200 32 32 0.0000000 6.2831853}%
% CIRCLE 2 0 3 0 Black White  
% 4 600 1200 630 1210 630 1210 630 1210
% 
\special{pn 8}%
\special{ar 600 1200 32 32 0.0000000 6.2831853}%
% LINE 2 0 3 0 Black White  
% 2 1790 1230 1400 1780
% 
\special{pn 8}%
\special{pa 1790 1230}%
\special{pa 1400 1780}%
\special{fp}%
% LINE 2 0 3 0 Black White  
% 2 1390 1230 1000 1800
% 
\special{pn 8}%
\special{pa 1390 1230}%
\special{pa 1000 1800}%
\special{fp}%
% LINE 2 0 3 0 Black White  
% 2 640 1230 1000 1800
% 
\special{pn 8}%
\special{pa 640 1230}%
\special{pa 1000 1800}%
\special{fp}%
% LINE 2 0 3 0 Black White  
% 2 1030 1230 1400 1780
% 
\special{pn 8}%
\special{pa 1030 1230}%
\special{pa 1400 1780}%
\special{fp}%
% STR 2 0 3 0 Black White  
% 4 1780 1044 1780 1110 2 0 0 0
% $\chi_3$
\put(17.8000,-11.1000){\makebox(0,0)[lb]{$\chi_3$}}%
\end{picture}}%
\end{center}
\bigskip
$\mathcal{F}(\widehat{G} \cup \widehat{G_0}) 
= \{(Ch(\chi_i), \circ), (Ch(\tau_j), \bullet): \chi_i \in \widehat{G}, \tau_j \in \widehat{G_0}\}$. 
Put $\gamma_i = (Ch(\chi_i), \circ)$ and $\rho_j = (Ch(\tau_j), \bullet)$.  
Then the structure equations are
\begin{align*}
&\gamma_1 \gamma_1 = \gamma_2,~~~\gamma_2 \gamma_2 = \gamma_0,~~~
\gamma_3 \gamma_3 = \gamma_1,~~~
\gamma_1 \gamma_2 = \gamma_3,~~~
\gamma_1 \gamma_3 = \gamma_0,~~~
\gamma_2 \gamma_3 = \gamma_1,\\
&\rho_0 \rho_0 = \rho_1 \rho_1 = \gamma_0 + \gamma_2,~~~
\rho_0 \rho_1 = \rho_1 \rho_0 = \gamma_1 + \gamma_3,~~~
\gamma_0 \rho_0 = \rho_0,~~~\gamma_1 \rho_0 = \rho_1,\\
&\gamma_2 \rho_0 = \rho_0,~~~\gamma_3 \rho_0 = \rho_1,~~~
\gamma_0 \rho_1 = \rho_1,~~~\gamma_1 \rho_1 = \rho_0,~~~
\gamma_2 \rho_1 = \rho_1,~~~\gamma_3 \rho_1 = \rho_0. 
\end{align*}
As in the case that example 4.1, we obtain a hypergroup isomorphic to 
the hypergroup of [HKTY1, example 4.4]. We write that only dimension function of 
$\mathcal{F}(\widehat{G} \cup \widehat{G_0})$. 
\begin{align*}
d(\gamma_0) = d(\gamma_1) = d(\gamma_2) = d(\gamma_3) =1,~~~
d(\rho_0) = d(\rho_1) = \sqrt{2}. 
\end{align*}

\bigskip
\noindent
{\bf 4.5 } \ The case that 
$G = \mathbb{Z}_2 \times \mathbb{Z}_2 = \{(e,e),(e,g), (g,e), (g,g)\}$ $(g^2 = e)$ and 
$G_0 = \mathbb{Z}_2 $.
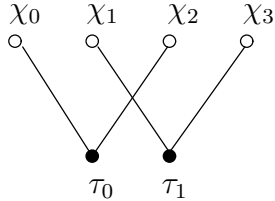
\begin{center}
%WinTpicVersion4.32a
{\unitlength 0.1in%
\begin{picture}(12.6400,9.1000)(5.6800,-18.9000)%
% STR 2 0 3 0 Black White  
% 4 570 1043 570 1110 2 0 0 0
% $\chi_0$
\put(5.7000,-11.1000){\makebox(0,0)[lb]{$\chi_0$}}%
% STR 2 0 3 0 Black White  
% 4 970 1044 970 1110 2 0 0 0
% $\chi_1$
\put(9.7000,-11.1000){\makebox(0,0)[lb]{$\chi_1$}}%
% STR 2 0 3 0 Black White  
% 4 980 1954 980 2020 2 0 0 0
% $\tau_0$
\put(9.8000,-20.2000){\makebox(0,0)[lb]{$\tau_0$}}%
% STR 2 0 3 0 Black White  
% 4 1360 1953 1360 2020 2 0 0 0
% $\tau_1$
\put(13.6000,-20.2000){\makebox(0,0)[lb]{$\tau_1$}}%
% STR 2 0 3 0 Black White  
% 4 1390 1044 1390 1110 2 0 0 0
% $\chi_2$
\put(13.9000,-11.1000){\makebox(0,0)[lb]{$\chi_2$}}%
% CIRCLE 2 0 3 0 Black White  
% 4 1000 1200 1030 1210 1030 1210 1030 1210
% 
\special{pn 8}%
\special{ar 1000 1200 32 32 0.0000000 6.2831853}%
% CIRCLE 2 0 0 0 Black Black  
% 4 1400 1800 1430 1810 1430 1810 1430 1810
% 
\special{sh 1.000}%
\special{ia 1400 1800 32 32 0.0000000 6.2831853}%
\special{pn 8}%
\special{ar 1400 1800 32 32 0.0000000 6.2831853}%
% CIRCLE 2 0 3 0 Black White  
% 4 1400 1200 1430 1210 1430 1210 1430 1210
% 
\special{pn 8}%
\special{ar 1400 1200 32 32 0.0000000 6.2831853}%
% CIRCLE 2 0 0 0 Black Black  
% 4 1000 1800 1030 1810 1030 1810 1030 1810
% 
\special{sh 1.000}%
\special{ia 1000 1800 32 32 0.0000000 6.2831853}%
\special{pn 8}%
\special{ar 1000 1800 32 32 0.0000000 6.2831853}%
% CIRCLE 2 0 3 0 Black Black  
% 4 1800 1200 1830 1210 1830 1210 1830 1210
% 
\special{pn 8}%
\special{ar 1800 1200 32 32 0.0000000 6.2831853}%
% CIRCLE 2 0 3 0 Black White  
% 4 600 1200 630 1210 630 1210 630 1210
% 
\special{pn 8}%
\special{ar 600 1200 32 32 0.0000000 6.2831853}%
% LINE 2 0 3 0 Black White  
% 2 1790 1230 1400 1780
% 
\special{pn 8}%
\special{pa 1790 1230}%
\special{pa 1400 1780}%
\special{fp}%
% LINE 2 0 3 0 Black White  
% 2 1390 1230 1000 1800
% 
\special{pn 8}%
\special{pa 1390 1230}%
\special{pa 1000 1800}%
\special{fp}%
% LINE 2 0 3 0 Black White  
% 2 640 1230 1000 1800
% 
\special{pn 8}%
\special{pa 640 1230}%
\special{pa 1000 1800}%
\special{fp}%
% LINE 2 0 3 0 Black White  
% 2 1030 1230 1400 1780
% 
\special{pn 8}%
\special{pa 1030 1230}%
\special{pa 1400 1780}%
\special{fp}%
% STR 2 0 3 0 Black White  
% 4 1780 1044 1780 1110 2 0 0 0
% $\chi_3$
\put(17.8000,-11.1000){\makebox(0,0)[lb]{$\chi_3$}}%
\end{picture}}%
\end{center}
\bigskip
$\mathcal{F}(\widehat{G} \cup \widehat{G_0}) 
= \{(Ch(\chi_i), \circ), (Ch(\tau_j), \bullet): \chi_i \in \widehat{G}, \tau_j \in \widehat{G_0}\}$. 
Put $\gamma_i = (Ch(\chi_i), \circ)$ and $\rho_j = (Ch(\tau_j), \bullet)$.  
Then the structure equations are
\begin{align*}
&\gamma_1 \gamma_1 = \gamma_0,~~~\gamma_2 \gamma_2 = \gamma_0,~~~
\gamma_3 \gamma_3 = \gamma_0,~~~
\gamma_1 \gamma_2 = \gamma_3,~~~
\gamma_1 \gamma_3 = \gamma_2,~~~
\gamma_2 \gamma_3 = \gamma_1,\\
&\rho_0 \rho_0 = \rho_1 \rho_1 = \gamma_0 + \gamma_2,~~~
\rho_0 \rho_1 = \rho_1 \rho_0 = \gamma_1 + \gamma_3,~~~
\gamma_0 \rho_0 = \rho_0,~~~\gamma_1 \rho_0 = \rho_1,\\
&\gamma_2 \rho_0 = \rho_0,~~~\gamma_3 \rho_0 = \rho_1,~~~
\gamma_0 \rho_1 = \rho_1,~~~\gamma_1 \rho_1 = \rho_0,~~~
\gamma_2 \rho_1 = \rho_1,~~~\gamma_3 \rho_1 = \rho_0. 
\end{align*}

\bigskip
\noindent
{\bf Remark} \ We note that Frobenius diagrams of 4.4 and 4.5 are same but their  
fusion rule algebra structures are different. 

\medskip
As in the case that example 4.1, we obtain a hypergroup isomorphic to 
the hypergroup of [HKTY1, example 4.5]. We write that only dimension function of 
$\mathcal{F}(\widehat{G} \cup \widehat{G_0})$. 
\begin{align*}
d(\gamma_0) = d(\gamma_1) = d(\gamma_2) = d(\gamma_3) =1,~~~
d(\rho_0) = d(\rho_1) = \sqrt{2}. 
\end{align*}

\bigskip
\noindent
{\bf 4.6 } \ The case that 
$G = S_3 = \mathbb{Z}_3 \rtimes_\alpha \mathbb{Z}_2$ and $G_0 = \mathbb{Z}_3$.
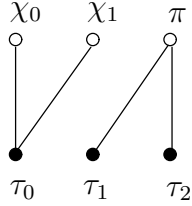
\begin{center}
%WinTpicVersion4.30a
{\unitlength 0.1in%
\begin{picture}(  8.6400,  9.3000)(  5.6800,-19.1000)%
% STR 2 0 3 0 Black White
% 4 570 1043 570 1110 2 0 0 0
% $\chi_0$
\put(5.7000,-11.1000){\makebox(0,0)[lb]{$\chi_0$}}%
% STR 2 0 3 0 Black White
% 4 970 1044 970 1110 2 0 0 0
% $\chi_1$
\put(9.7000,-11.1000){\makebox(0,0)[lb]{$\chi_1$}}%
% STR 2 0 3 0 Black White
% 4 570 1964 570 2030 2 0 0 0
% $\tau_0$
\put(5.7000,-20.3000){\makebox(0,0)[lb]{$\tau_0$}}%
% STR 2 0 3 0 Black White
% 4 950 1963 950 2030 2 0 0 0
% $\tau_1$
\put(9.5000,-20.3000){\makebox(0,0)[lb]{$\tau_1$}}%
% STR 2 0 3 0 Black White
% 4 1390 1044 1390 1110 2 0 0 0
% $\pi$
\put(13.9000,-11.1000){\makebox(0,0)[lb]{$\pi$}}%
% CIRCLE 2 0 3 0 Black White
% 4 1000 1200 1030 1210 1030 1210 1030 1210
% 
\special{pn 8}%
\special{ar 1000 1200 32 32  0.0000000  6.2831853}%
% CIRCLE 2 0 0 0 Black Black
% 4 600 1800 630 1810 630 1810 630 1810
% 
\special{sh 1.000}%
\special{ia 600 1800 32 32  0.0000000  6.2831853}%
\special{pn 8}%
\special{ar 600 1800 32 32  0.0000000  6.2831853}%
% CIRCLE 2 0 3 0 Black White
% 4 1400 1200 1430 1210 1430 1210 1430 1210
% 
\special{pn 8}%
\special{ar 1400 1200 32 32  0.0000000  6.2831853}%
% CIRCLE 2 0 0 0 Black Black
% 4 1000 1800 1030 1810 1030 1810 1030 1810
% 
\special{sh 1.000}%
\special{ia 1000 1800 32 32  0.0000000  6.2831853}%
\special{pn 8}%
\special{ar 1000 1800 32 32  0.0000000  6.2831853}%
% CIRCLE 2 0 0 0 Black Black
% 4 1400 1800 1430 1810 1430 1810 1430 1810
% 
\special{sh 1.000}%
\special{ia 1400 1800 32 32  0.0000000  6.2831853}%
\special{pn 8}%
\special{ar 1400 1800 32 32  0.0000000  6.2831853}%
% CIRCLE 2 0 3 0 Black White
% 4 600 1200 630 1210 630 1210 630 1210
% 
\special{pn 8}%
\special{ar 600 1200 32 32  0.0000000  6.2831853}%
% LINE 2 0 3 0 Black White
% 4 600 1240 600 1790 600 1790 600 1770
% 
\special{pn 8}%
\special{pa 600 1240}%
\special{pa 600 1790}%
\special{fp}%
\special{pa 600 1790}%
\special{pa 600 1770}%
\special{fp}%
% LINE 2 0 3 0 Black White
% 2 990 1240 600 1790
% 
\special{pn 8}%
\special{pa 990 1240}%
\special{pa 600 1790}%
\special{fp}%
% LINE 2 0 3 0 Black White
% 2 1390 1230 1000 1800
% 
\special{pn 8}%
\special{pa 1390 1230}%
\special{pa 1000 1800}%
\special{fp}%
% LINE 2 0 3 0 Black White
% 2 1410 1240 1410 1780
% 
\special{pn 8}%
\special{pa 1410 1240}%
\special{pa 1410 1780}%
\special{fp}%
% STR 2 0 3 0 Black White
% 4 1380 1973 1380 2040 2 0 0 0
% $\tau_2$
\put(13.8000,-20.4000){\makebox(0,0)[lb]{$\tau_2$}}%
\end{picture}}%
\end{center}
\bigskip
$\mathcal{F}(\widehat{G} \cup \widehat{G_0})$ is not a 
fusion rule algebra by Corollary 3.17.

\bigskip
\noindent
{\bf 4.7 } \ The case that $G$ is the dihedral group 
$D_4 = \mathbb{Z}_4 \rtimes_\alpha \mathbb{Z}_2$ and $G_0 = \mathbb{Z}_2$.
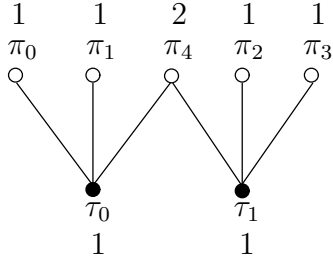
\begin{center}
%WinTpicVersion4.32a
{\unitlength 0.1in%
\begin{picture}(16.1000,12.0000)(19.8000,-19.1000)%
% STR 2 0 3 0 Black White  
% 4 1980 951 1980 1021 2 0 0 0
% $\pi_0$
\put(19.8000,-10.2100){\makebox(0,0)[lb]{$\pi_0$}}%
% STR 2 0 3 0 Black White  
% 4 2388 951 2388 1021 2 0 0 0
% $\pi_1$
\put(23.8800,-10.2100){\makebox(0,0)[lb]{$\pi_1$}}%
% STR 2 0 3 0 Black White  
% 4 2381 1794 2381 1865 2 0 0 0
% $\tau_0$
\put(23.8100,-18.6500){\makebox(0,0)[lb]{$\tau_0$}}%
% CIRCLE 2 0 3 0 Black White  
% 4 2423 1106 2425 1070 2425 1070 2425 1070
% 
\special{pn 8}%
\special{ar 2423 1106 36 36 0.0000000 6.2831853}%
% CIRCLE 2 0 3 0 Black White  
% 4 2022 1106 2024 1070 2024 1070 2024 1070
% 
\special{pn 8}%
\special{ar 2022 1106 36 36 0.0000000 6.2831853}%
% CIRCLE 2 0 0 0 Black Black  
% 4 2423 1703 2425 1667 2425 1667 2425 1667
% 
\special{sh 1.000}%
\special{ia 2423 1703 36 36 0.0000000 6.2831853}%
\special{pn 8}%
\special{ar 2423 1703 36 36 0.0000000 6.2831853}%
% CIRCLE 2 0 3 0 Black White  
% 4 2830 1110 2832 1074 2832 1074 2832 1074
% 
\special{pn 8}%
\special{ar 2830 1110 36 36 0.0000000 6.2831853}%
% LINE 2 0 3 0 Black White  
% 6 2423 1148 2423 1682 2823 1148 2423 1682 2029 1148 2423 1682
% 
\special{pn 8}%
\special{pa 2423 1148}%
\special{pa 2423 1682}%
\special{fp}%
\special{pa 2823 1148}%
\special{pa 2423 1682}%
\special{fp}%
\special{pa 2029 1148}%
\special{pa 2423 1682}%
\special{fp}%
% STR 2 0 3 0 Black White  
% 4 3154 951 3154 1021 2 0 0 0
% $\pi_2$
\put(31.5400,-10.2100){\makebox(0,0)[lb]{$\pi_2$}}%
% CIRCLE 2 0 3 0 Black White  
% 4 3196 1106 3198 1070 3198 1070 3198 1070
% 
\special{pn 8}%
\special{ar 3196 1106 36 36 0.0000000 6.2831853}%
% CIRCLE 2 0 3 0 Black White  
% 4 3554 1106 3556 1070 3556 1070 3556 1070
% 
\special{pn 8}%
\special{ar 3554 1106 36 36 0.0000000 6.2831853}%
% CIRCLE 2 0 0 0 Black Black  
% 4 3196 1710 3198 1674 3198 1674 3198 1674
% 
\special{sh 1.000}%
\special{ia 3196 1710 36 36 0.0000000 6.2831853}%
\special{pn 8}%
\special{ar 3196 1710 36 36 0.0000000 6.2831853}%
% LINE 2 0 3 0 Black Black  
% 6 2837 1148 3189 1675 3196 1148 3196 1675 3561 1148 3196 1675
% 
\special{pn 8}%
\special{pa 2837 1148}%
\special{pa 3189 1675}%
\special{fp}%
\special{pa 3196 1148}%
\special{pa 3196 1675}%
\special{fp}%
\special{pa 3561 1148}%
\special{pa 3196 1675}%
\special{fp}%
% STR 2 0 3 0 Black White  
% 4 2788 951 2788 1021 2 0 0 0
% $\pi_4$
\put(27.8800,-10.2100){\makebox(0,0)[lb]{$\pi_4$}}%
% STR 2 0 3 0 Black White  
% 4 3512 951 3512 1021 2 0 0 0
% $\pi_3$
\put(35.1200,-10.2100){\makebox(0,0)[lb]{$\pi_3$}}%
% STR 2 0 3 0 Black White  
% 4 3154 1794 3154 1865 2 0 0 0
% $\tau_1$
\put(31.5400,-18.6500){\makebox(0,0)[lb]{$\tau_1$}}%
% STR 2 0 3 0 Black White  
% 4 2000 740 2000 840 2 0 0 0
% 1
\put(20.0000,-8.4000){\makebox(0,0)[lb]{1}}%
% STR 2 0 3 0 Black White  
% 4 2420 740 2420 840 2 0 0 0
% 1
\put(24.2000,-8.4000){\makebox(0,0)[lb]{1}}%
% STR 2 0 3 0 Black White  
% 4 2830 740 2830 840 2 0 0 0
% 2
\put(28.3000,-8.4000){\makebox(0,0)[lb]{2}}%
% STR 2 0 3 0 Black White  
% 4 3190 740 3190 840 2 0 0 0
% 1
\put(31.9000,-8.4000){\makebox(0,0)[lb]{1}}%
% STR 2 0 3 0 Black White  
% 4 3550 740 3550 840 2 0 0 0
% 1
\put(35.5000,-8.4000){\makebox(0,0)[lb]{1}}%
% STR 2 0 3 0 Black White  
% 4 2410 1940 2410 2040 2 0 0 0
% 1
\put(24.1000,-20.4000){\makebox(0,0)[lb]{1}}%
% STR 2 0 3 0 Black White  
% 4 3180 1940 3180 2040 2 0 0 0
% 1
\put(31.8000,-20.4000){\makebox(0,0)[lb]{1}}%
\end{picture}}%
\end{center}
\bigskip
$\mathcal{F}(\widehat{G} \cup \widehat{G_0}) 
= \{(Ch(\pi_i), \circ), (Ch(\tau_j), \bullet): \pi_i \in \widehat{G}, \tau_j \in \widehat{G_0}\}$. 
Put $\gamma_i = (Ch(\pi_i), \circ)$ and $\rho_j = (Ch(\tau_j), \bullet)$.  
Then the structure equations are
\begin{align*}
&\gamma_1 \gamma_1 = \gamma_0,~~~ 
\gamma_2 \gamma_2 
= \gamma_0 + \gamma_1 + \gamma_3 + \gamma_4,~~~
\gamma_3 \gamma_3 = \gamma_0,~~~
\gamma_4 \gamma_4 = \gamma_0,\\
&\gamma_1 \gamma_2 = \gamma_2,~~~
\gamma_1 \gamma_3 = \gamma_4,~~~
\gamma_1 \gamma_4 = \gamma_3,~~~
\gamma_2 \gamma_3 = \gamma_2,~~~
\gamma_2 \gamma_4 = \gamma_2,~~~
\gamma_3 \gamma_4 = \gamma_1,\\
&\rho_0 \rho_0 = \rho_1 \rho_1 = \gamma_0 + \gamma_1 + \gamma_2,~~~
\rho_0 \rho_1 = \gamma_2 + \gamma_3 + \gamma_4,\\
&\gamma_1 \rho_0 = \rho_0,~~~
\gamma_2 \rho_0 = \rho_0 + \rho_1,
\gamma_3 \rho_0 = \rho_1,~~~
\gamma_4 \rho_0 = \rho_1,\\
&\gamma_1 \rho_1 = \rho_1,~~~
\gamma_2 \rho_1 = \rho_0 + \rho_1,~~~
\gamma_3 \rho_1 = \rho_0,~~~
\gamma_4 \rho_1 = \rho_0.
\end{align*}  
As in the case that example 4.1, we obtain a hypergroup isomorphic to 
the hypergroup of [HKTY1, example 4.7]. We write that only dimension function of 
$\mathcal{F}(\widehat{G} \cup \widehat{G_0})$. 
\begin{align*}
d(\gamma_0) = d(\gamma_1) = d(\gamma_2) = d(\gamma_3) =1,~~~
d(\gamma_4) = d(\rho_0) = d(\rho_1) = 2. 
\end{align*}

\bigskip
\noindent
{\bf 4.8 } \ The case that $G$ is the alternating group 
$A_4 = (\mathbb{Z}_2 \times \mathbb{Z}_2) \rtimes_\alpha \mathbb{Z}_3$ of 
degree 4 and  $G_0 = \mathbb{Z}_3$.
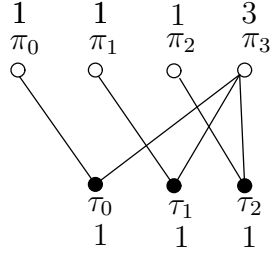
\begin{center}
%WinTpicVersion4.32a
{\unitlength 0.1in%
\begin{picture}(12.5200,11.7000)(19.8000,-19.0000)%
% STR 2 0 3 0 Black White  
% 4 1980 951 1980 1021 2 0 0 0
% $\pi_0$
\put(19.8000,-10.2100){\makebox(0,0)[lb]{$\pi_0$}}%
% STR 2 0 3 0 Black White  
% 4 2388 951 2388 1021 2 0 0 0
% $\pi_1$
\put(23.8800,-10.2100){\makebox(0,0)[lb]{$\pi_1$}}%
% STR 2 0 3 0 Black White  
% 4 2381 1794 2381 1865 2 0 0 0
% $\tau_0$
\put(23.8100,-18.6500){\makebox(0,0)[lb]{$\tau_0$}}%
% CIRCLE 2 0 3 0 Black White  
% 4 2423 1106 2425 1070 2425 1070 2425 1070
% 
\special{pn 8}%
\special{ar 2423 1106 36 36 0.0000000 6.2831853}%
% CIRCLE 2 0 3 0 Black White  
% 4 2022 1106 2024 1070 2024 1070 2024 1070
% 
\special{pn 8}%
\special{ar 2022 1106 36 36 0.0000000 6.2831853}%
% CIRCLE 2 0 0 0 Black Black  
% 4 2423 1703 2425 1667 2425 1667 2425 1667
% 
\special{sh 1.000}%
\special{ia 2423 1703 36 36 0.0000000 6.2831853}%
\special{pn 8}%
\special{ar 2423 1703 36 36 0.0000000 6.2831853}%
% CIRCLE 2 0 3 0 Black White  
% 4 2830 1110 2832 1074 2832 1074 2832 1074
% 
\special{pn 8}%
\special{ar 2830 1110 36 36 0.0000000 6.2831853}%
% STR 2 0 3 0 Black White  
% 4 2790 940 2790 1010 2 0 0 0
% $\pi_2$
\put(27.9000,-10.1000){\makebox(0,0)[lb]{$\pi_2$}}%
% CIRCLE 2 0 3 0 Black White  
% 4 3196 1106 3198 1070 3198 1070 3198 1070
% 
\special{pn 8}%
\special{ar 3196 1106 36 36 0.0000000 6.2831853}%
% CIRCLE 2 0 0 0 Black Black  
% 4 2830 1710 2832 1674 2832 1674 2832 1674
% 
\special{sh 1.000}%
\special{ia 2830 1710 36 36 0.0000000 6.2831853}%
\special{pn 8}%
\special{ar 2830 1710 36 36 0.0000000 6.2831853}%
% CIRCLE 2 0 0 0 Black Black  
% 4 3196 1710 3198 1674 3198 1674 3198 1674
% 
\special{sh 1.000}%
\special{ia 3196 1710 36 36 0.0000000 6.2831853}%
\special{pn 8}%
\special{ar 3196 1710 36 36 0.0000000 6.2831853}%
% STR 2 0 3 0 Black White  
% 4 3180 940 3180 1010 2 0 0 0
% $\pi_3$
\put(31.8000,-10.1000){\makebox(0,0)[lb]{$\pi_3$}}%
% STR 2 0 3 0 Black White  
% 4 3154 1794 3154 1865 2 0 0 0
% $\tau_2$
\put(31.5400,-18.6500){\makebox(0,0)[lb]{$\tau_2$}}%
% STR 2 0 3 0 Black White  
% 4 2800 1809 2800 1880 2 0 0 0
% $\tau_1$
\put(28.0000,-18.8000){\makebox(0,0)[lb]{$\tau_1$}}%
% LINE 2 0 3 0 Black White  
% 2 2030 1140 2423 1703
% 
\special{pn 8}%
\special{pa 2030 1140}%
\special{pa 2423 1703}%
\special{fp}%
% LINE 2 0 3 0 Black White  
% 2 2440 1140 2830 1710
% 
\special{pn 8}%
\special{pa 2440 1140}%
\special{pa 2830 1710}%
\special{fp}%
% LINE 2 0 3 0 Black White  
% 2 2860 1150 3196 1710
% 
\special{pn 8}%
\special{pa 2860 1150}%
\special{pa 3196 1710}%
\special{fp}%
% LINE 2 0 3 0 Black White  
% 2 3170 1130 2423 1703
% 
\special{pn 8}%
\special{pa 3170 1130}%
\special{pa 2423 1703}%
\special{fp}%
% LINE 2 0 3 0 Black White  
% 2 3170 1130 2830 1710
% 
\special{pn 8}%
\special{pa 3170 1130}%
\special{pa 2830 1710}%
\special{fp}%
% LINE 2 0 3 0 Black White  
% 2 3170 1130 3196 1710
% 
\special{pn 8}%
\special{pa 3170 1130}%
\special{pa 3196 1710}%
\special{fp}%
% STR 2 0 3 0 Black White  
% 4 1990 760 1990 860 2 0 0 0
% 1
\put(19.9000,-8.6000){\makebox(0,0)[lb]{1}}%
% STR 2 0 3 0 Black White  
% 4 1990 760 1990 860 2 0 0 0
% 1
\put(19.9000,-8.6000){\makebox(0,0)[lb]{1}}%
% STR 2 0 3 0 Black White  
% 4 2410 760 2410 860 2 0 0 0
% 1
\put(24.1000,-8.6000){\makebox(0,0)[lb]{1}}%
% STR 2 0 3 0 Black White  
% 4 2410 760 2410 860 2 0 0 0
% 1
\put(24.1000,-8.6000){\makebox(0,0)[lb]{1}}%
% STR 2 0 3 0 Black White  
% 4 2810 770 2810 870 2 0 0 0
% 1
\put(28.1000,-8.7000){\makebox(0,0)[lb]{1}}%
% STR 2 0 3 0 Black White  
% 4 2410 1920 2410 2020 2 0 0 0
% 1
\put(24.1000,-20.2000){\makebox(0,0)[lb]{1}}%
% STR 2 0 3 0 Black White  
% 4 2830 1930 2830 2030 2 0 0 0
% 1
\put(28.3000,-20.3000){\makebox(0,0)[lb]{1}}%
% STR 2 0 3 0 Black White  
% 4 3180 1930 3180 2030 2 0 0 0
% 1
\put(31.8000,-20.3000){\makebox(0,0)[lb]{1}}%
% STR 2 0 3 0 Black White  
% 4 3180 760 3180 860 2 0 0 0
% 3
\put(31.8000,-8.6000){\makebox(0,0)[lb]{3}}%
\end{picture}}%
\end{center}
\bigskip
$\mathcal{F}(\widehat{G} \cup \widehat{G_0}) 
= \{(Ch(\pi_i), \circ), (Ch(\tau_j), \bullet): \pi_i \in \widehat{G}, \tau_j \in \widehat{G_0}\}$. 
Put $\gamma_i = (Ch(\pi_i), \circ)$ and $\rho_j = (Ch(\tau_j), \bullet)$.
Then the structure equations are
\begin{align*}
&\gamma_1 \gamma_1 = \gamma_2,~~~
\gamma_2 \gamma_2 = \gamma_1,~~~
\gamma_3 \gamma_3 
= \gamma_0 + \gamma_1 + \gamma_2 + 2\gamma_3,~~~
\gamma_1 \gamma_2 = \gamma_0,~~~
\gamma_1 \gamma_3 = \gamma_3,\\
&\gamma_2 \gamma_3 = \gamma_3,~~~
\rho_0 \rho_0 = \rho_1 \rho_2 = \gamma_0 +  \gamma_3,~~~
\rho_0 \rho_1 = \gamma_1 + \gamma_3,~~~
\rho_0 \rho_2 = \gamma_2 + \gamma_3,\\
&\gamma_1 \rho_0 = \rho_1,~~~
\gamma_2 \rho_0 = \rho_2,~~~
\gamma_1 \rho_1 = \rho_2,~~~
\gamma_2 \rho_1 = \rho_0, \\
&\gamma_3 \rho_0 = \gamma_3 \rho_1 = \gamma_3 \rho_2 = \rho_0 + \rho_1 + \rho_2. 
\end{align*}  
As in the case that example 4.1, we obtain a hypergroup isomorphic to 
the hypergroup of [HKTY1, example 4.8]. We write that only dimension function of 
$\mathcal{F}(\widehat{G} \cup \widehat{G_0})$. 
\begin{align*}
d(\gamma_0) = d(\gamma_1) = d(\gamma_2) = 1,~~~d(\gamma_3) =3,~~~
d(\rho_0) = d(\rho_1) = d(\rho_2) = 2. 
\end{align*}

\bigskip
\noindent
{\bf 4.9 } \ The case that $G$ is the symmetric group 
$S_4 = A_4 \rtimes_\alpha \mathbb{Z}_2$ of degree 4 and $G_0 = \mathbb{Z}_2$.
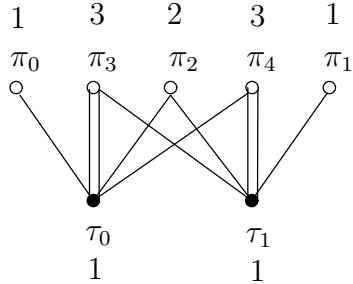
\begin{center}
%WinTpicVersion4.32a
{\unitlength 0.1in%
\begin{picture}(16.8600,13.4000)(5.6800,-20.8000)%
% STR 2 0 3 0 Black White  
% 4 570 1043 570 1110 2 0 0 0
% $\pi_0$
\put(5.7000,-11.1000){\makebox(0,0)[lb]{$\pi_0$}}%
% STR 2 0 3 0 Black White  
% 4 970 1044 970 1110 2 0 0 0
% $\pi_3$
\put(9.7000,-11.1000){\makebox(0,0)[lb]{$\pi_3$}}%
% STR 2 0 3 0 Black White  
% 4 960 1954 960 2020 2 0 0 0
% $\tau_0$
\put(9.6000,-20.2000){\makebox(0,0)[lb]{$\tau_0$}}%
% STR 2 0 3 0 Black White  
% 4 1790 1963 1790 2030 2 0 0 0
% $\tau_1$
\put(17.9000,-20.3000){\makebox(0,0)[lb]{$\tau_1$}}%
% STR 2 0 3 0 Black White  
% 4 1390 1044 1390 1110 2 0 0 0
% $\pi_2$
\put(13.9000,-11.1000){\makebox(0,0)[lb]{$\pi_2$}}%
% CIRCLE 2 0 3 0 Black White  
% 4 1000 1200 1030 1210 1030 1210 1030 1210
% 
\special{pn 8}%
\special{ar 1000 1200 32 32 0.0000000 6.2831853}%
% CIRCLE 2 0 0 0 Black Black  
% 4 1000 1790 1030 1800 1030 1800 1030 1800
% 
\special{sh 1.000}%
\special{ia 1000 1790 32 32 0.0000000 6.2831853}%
\special{pn 8}%
\special{ar 1000 1790 32 32 0.0000000 6.2831853}%
% CIRCLE 2 0 3 0 Black White  
% 4 1400 1200 1430 1210 1430 1210 1430 1210
% 
\special{pn 8}%
\special{ar 1400 1200 32 32 0.0000000 6.2831853}%
% CIRCLE 2 0 0 0 Black Black  
% 4 1820 1790 1850 1800 1850 1800 1850 1800
% 
\special{sh 1.000}%
\special{ia 1820 1790 32 32 0.0000000 6.2831853}%
\special{pn 8}%
\special{ar 1820 1790 32 32 0.0000000 6.2831853}%
% CIRCLE 2 0 3 0 Black White  
% 4 600 1200 630 1210 630 1210 630 1210
% 
\special{pn 8}%
\special{ar 600 1200 32 32 0.0000000 6.2831853}%
% STR 2 0 3 0 Black White  
% 4 1792 1043 1792 1110 2 0 0 0
% $\pi_4$
\put(17.9200,-11.1000){\makebox(0,0)[lb]{$\pi_4$}}%
% STR 2 0 3 0 Black White  
% 4 2192 1044 2192 1110 2 0 0 0
% $\pi_1$
\put(21.9200,-11.1000){\makebox(0,0)[lb]{$\pi_1$}}%
% CIRCLE 2 0 3 0 Black White  
% 4 2222 1200 2252 1210 2252 1210 2252 1210
% 
\special{pn 8}%
\special{ar 2222 1200 32 32 0.0000000 6.2831853}%
% CIRCLE 2 0 3 0 Black White  
% 4 1822 1200 1852 1210 1852 1210 1852 1210
% 
\special{pn 8}%
\special{ar 1822 1200 32 32 0.0000000 6.2831853}%
% LINE 2 0 3 0 Black White  
% 2 1030 1210 1030 1800
% 
\special{pn 8}%
\special{pa 1030 1210}%
\special{pa 1030 1800}%
\special{fp}%
% LINE 2 0 3 0 Black White  
% 2 1800 1200 1800 1800
% 
\special{pn 8}%
\special{pa 1800 1200}%
\special{pa 1800 1800}%
\special{fp}%
% LINE 2 0 3 0 Black White  
% 2 1852 1210 1850 1800
% 
\special{pn 8}%
\special{pa 1852 1210}%
\special{pa 1850 1800}%
\special{fp}%
% LINE 2 0 3 0 Black White  
% 2 980 1210 980 1780
% 
\special{pn 8}%
\special{pa 980 1210}%
\special{pa 980 1780}%
\special{fp}%
% LINE 2 0 3 0 Black White  
% 2 1400 1240 1000 1790
% 
\special{pn 8}%
\special{pa 1400 1240}%
\special{pa 1000 1790}%
\special{fp}%
% LINE 2 0 3 0 Black White  
% 2 1800 1230 1000 1790
% 
\special{pn 8}%
\special{pa 1800 1230}%
\special{pa 1000 1790}%
\special{fp}%
% LINE 2 0 3 0 Black White  
% 2 1030 1210 1820 1790
% 
\special{pn 8}%
\special{pa 1030 1210}%
\special{pa 1820 1790}%
\special{fp}%
% LINE 2 0 3 0 Black White  
% 2 1400 1240 1820 1790
% 
\special{pn 8}%
\special{pa 1400 1240}%
\special{pa 1820 1790}%
\special{fp}%
% LINE 2 0 3 0 Black White  
% 2 2210 1230 1820 1790
% 
\special{pn 8}%
\special{pa 2210 1230}%
\special{pa 1820 1790}%
\special{fp}%
% LINE 2 0 3 0 Black White  
% 2 620 1240 1000 1790
% 
\special{pn 8}%
\special{pa 620 1240}%
\special{pa 1000 1790}%
\special{fp}%
% STR 2 0 3 0 Black White  
% 4 570 790 570 890 2 0 0 0
% 1
\put(5.7000,-8.9000){\makebox(0,0)[lb]{1}}%
% STR 2 0 3 0 Black White  
% 4 980 770 980 870 2 0 0 0
% 3
\put(9.8000,-8.7000){\makebox(0,0)[lb]{3}}%
% STR 2 0 3 0 Black White  
% 4 1380 770 1380 870 2 0 0 0
% 2
\put(13.8000,-8.7000){\makebox(0,0)[lb]{2}}%
% STR 2 0 3 0 Black White  
% 4 1810 780 1810 880 2 0 0 0
% 3
\put(18.1000,-8.8000){\makebox(0,0)[lb]{3}}%
% STR 2 0 3 0 Black White  
% 4 2200 770 2200 870 2 0 0 0
% 1
\put(22.0000,-8.7000){\makebox(0,0)[lb]{1}}%
% STR 2 0 3 0 Black White  
% 4 970 2100 970 2200 2 0 0 0
% 1
\put(9.7000,-22.0000){\makebox(0,0)[lb]{1}}%
% STR 2 0 3 0 Black White  
% 4 1810 2110 1810 2210 2 0 0 0
% 1
\put(18.1000,-22.1000){\makebox(0,0)[lb]{1}}%
\end{picture}}%
\end{center}
\bigskip
$\mathcal{F}(\widehat{G} \cup \widehat{G_0}) 
= \{(Ch(\pi_i), \circ), (Ch(\tau_j), \bullet): \pi_i \in \widehat{G}, \tau_j \in \widehat{G_0}\}$. 
Put $\gamma_i = (Ch(\pi_i), \circ)$ and $\rho_j = (Ch(\tau_j), \bullet)$.  
Then the structure equations are
\begin{align*}
&\gamma_1  \gamma_1 = \gamma_0,~~~
\gamma_2  \gamma_2 = \gamma_0 + \gamma_1 + \gamma_2,~~~
\gamma_3 \gamma_3 = \gamma_4  \gamma_4 
= \gamma_0 + \gamma_2 + \gamma_3 + \gamma_4,\\
&\gamma_1  \gamma_2 = \gamma_2,~~~\gamma_1  \gamma_3 = \gamma_4,~~~
\gamma_1 \gamma_4 = \gamma_3,~~~
\gamma_2 \gamma_3 = \gamma_2  \gamma_4 
= \gamma_3 + \gamma_4,\\
&\gamma_3 \gamma_4 
= \gamma_1 + \gamma_2 + \gamma_3 + \gamma_4,~~~
\rho_0 \rho_0 = \rho_1 \rho_1 = 
\gamma_0 + \gamma_2 + \gamma_3 + \gamma_4,\\
&\rho_0 \rho_1 = \rho_1 \rho_0 = 
\gamma_1 + \gamma_2 + \gamma_3 + \gamma_4,~~~
\gamma_0 \rho_0 = \rho_0,~~~\gamma_1 \rho_0 = \rho_1,\\
&\gamma_2 \rho_0 = \rho_0 + \rho_1,~~~  
\gamma_3 \rho_0 =  \rho_0 + \rho_1,~~~
\gamma_4 \rho_0 =  \rho_0 + \rho_1,~~~
\gamma_0 \rho_1 = \rho_1,\\
&\gamma_1 \rho_1 = \rho_0,~~~
\gamma_2 \rho_1 = \rho_0 + \rho_1,~~~  
\gamma_3 \rho_1 =  \rho_0 + \rho_1,~~~
\gamma_4 \rho_1 =  \rho_0 + \rho_1.          
\end{align*}
As in the case that example 4.1, we obtain a hypergroup isomorphic to 
the hypergroup of [HKTY1, example 4.9]. We write that only dimension function of 
$\mathcal{F}(\widehat{G} \cup \widehat{G_0})$. 
\begin{align*}
d(\gamma_0) = d(\gamma_1) = 1,~~~d(\gamma_2) = 2,~~~d(\gamma_3) = d(\gamma_4) = 3,~~~
d(\rho_0) = d(\rho_1) = 2 \sqrt{3}. 
\end{align*}

\bigskip
\noindent
{\bf 4.10 } \ The case that $G$ is the symmetric group 
$S_4$ of degree 4 and $G_0$ is the symmetric group $S_3$ of degree 3.
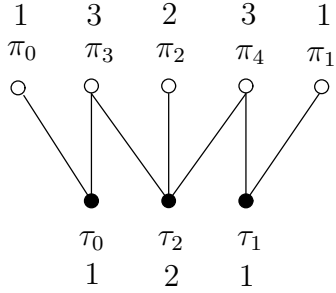
\begin{center}
%WinTpicVersion4.32a
{\unitlength 0.1in%
\begin{picture}(16.6400,13.7000)(17.7000,-19.1000)%
% STR 2 0 3 0 Black White  
% 4 1770 791 1770 860 2 0 0 0
% $\pi_0$
\put(17.7000,-8.6000){\makebox(0,0)[lb]{$\pi_0$}}%
% STR 2 0 3 0 Black White  
% 4 2950 811 2950 880 2 0 0 0
% $\pi_4$
\put(29.5000,-8.8000){\makebox(0,0)[lb]{$\pi_4$}}%
% STR 2 0 3 0 Black White  
% 4 2140 1780 2140 1850 2 0 0 0
% $\tau_0$
\put(21.4000,-18.5000){\makebox(0,0)[lb]{$\tau_0$}}%
% CIRCLE 2 0 3 0 Black White  
% 4 2606 996 2608 960 2608 960 2608 960
% 
\special{pn 8}%
\special{ar 2606 996 36 36 0.0000000 6.2831853}%
% CIRCLE 2 0 3 0 Black White  
% 4 3398 1000 3400 964 3400 964 3400 964
% 
\special{pn 8}%
\special{ar 3398 1000 36 36 0.0000000 6.2831853}%
% CIRCLE 2 0 0 0 Black Black  
% 4 2606 1596 2608 1561 2608 1561 2608 1561
% 
\special{sh 1.000}%
\special{ia 2606 1596 35 35 0.0000000 6.2831853}%
\special{pn 8}%
\special{ar 2606 1596 35 35 0.0000000 6.2831853}%
% CIRCLE 2 0 0 0 Black Black  
% 4 2206 1596 2208 1560 2208 1560 2208 1560
% 
\special{sh 1.000}%
\special{ia 2206 1596 36 36 0.0000000 6.2831853}%
\special{pn 8}%
\special{ar 2206 1596 36 36 0.0000000 6.2831853}%
% STR 2 0 3 0 Black White  
% 4 2550 1780 2550 1850 2 0 0 0
% $\tau_2$
\put(25.5000,-18.5000){\makebox(0,0)[lb]{$\tau_2$}}%
% STR 2 0 3 0 Black White  
% 4 2170 801 2170 870 2 0 0 0
% $\pi_3$
\put(21.7000,-8.7000){\makebox(0,0)[lb]{$\pi_3$}}%
% STR 2 0 3 0 Black White  
% 4 3320 821 3320 890 2 0 0 0
% $\pi_1$
\put(33.2000,-8.9000){\makebox(0,0)[lb]{$\pi_1$}}%
% CIRCLE 2 0 3 0 Black White  
% 4 2206 996 2208 960 2208 960 2208 960
% 
\special{pn 8}%
\special{ar 2206 996 36 36 0.0000000 6.2831853}%
% CIRCLE 2 0 3 0 Black White  
% 4 3006 996 3008 960 3008 960 3008 960
% 
\special{pn 8}%
\special{ar 3006 996 36 36 0.0000000 6.2831853}%
% CIRCLE 2 0 0 0 Black Black  
% 4 3006 1596 3008 1560 3008 1560 3008 1560
% 
\special{sh 1.000}%
\special{ia 3006 1596 36 36 0.0000000 6.2831853}%
\special{pn 8}%
\special{ar 3006 1596 36 36 0.0000000 6.2831853}%
% STR 2 0 3 0 Black White  
% 4 2960 1780 2960 1850 2 0 0 0
% $\tau_1$
\put(29.6000,-18.5000){\makebox(0,0)[lb]{$\tau_1$}}%
% STR 2 0 3 0 Black White  
% 4 2540 801 2540 870 2 0 0 0
% $\pi_2$
\put(25.4000,-8.7000){\makebox(0,0)[lb]{$\pi_2$}}%
% CIRCLE 2 0 3 0 Black White  
% 4 1826 1006 1828 970 1828 970 1828 970
% 
\special{pn 8}%
\special{ar 1826 1006 36 36 0.0000000 6.2831853}%
% LINE 2 0 3 0 Black Black  
% 2 1856 1046 2206 1596
% 
\special{pn 8}%
\special{pa 1856 1046}%
\special{pa 2206 1596}%
\special{fp}%
% LINE 2 0 3 0 Black Black  
% 2 2206 1036 2206 1596
% 
\special{pn 8}%
\special{pa 2206 1036}%
\special{pa 2206 1596}%
\special{fp}%
% LINE 2 0 3 0 Black Black  
% 2 2206 1036 2606 1596
% 
\special{pn 8}%
\special{pa 2206 1036}%
\special{pa 2606 1596}%
\special{fp}%
% LINE 2 0 3 0 Black Black  
% 2 2606 1036 2606 1596
% 
\special{pn 8}%
\special{pa 2606 1036}%
\special{pa 2606 1596}%
\special{fp}%
% LINE 2 0 3 0 Black Black  
% 2 3006 1036 3006 1596
% 
\special{pn 8}%
\special{pa 3006 1036}%
\special{pa 3006 1596}%
\special{fp}%
% LINE 2 0 3 0 Black Black  
% 2 3386 1036 3006 1596
% 
\special{pn 8}%
\special{pa 3386 1036}%
\special{pa 3006 1596}%
\special{fp}%
% LINE 2 0 3 0 Black Black  
% 2 3006 1036 2606 1596
% 
\special{pn 8}%
\special{pa 3006 1036}%
\special{pa 2606 1596}%
\special{fp}%
% STR 2 0 3 0 Black White  
% 4 1800 580 1800 680 2 0 0 0
% 1
\put(18.0000,-6.8000){\makebox(0,0)[lb]{1}}%
% STR 2 0 3 0 Black White  
% 4 2180 580 2180 680 2 0 0 0
% 3
\put(21.8000,-6.8000){\makebox(0,0)[lb]{3}}%
% STR 2 0 3 0 Black White  
% 4 2570 570 2570 670 2 0 0 0
% 2
\put(25.7000,-6.7000){\makebox(0,0)[lb]{2}}%
% STR 2 0 3 0 Black White  
% 4 2980 570 2980 670 2 0 0 0
% 3
\put(29.8000,-6.7000){\makebox(0,0)[lb]{3}}%
% STR 2 0 3 0 Black White  
% 4 3370 580 3370 680 2 0 0 0
% 1
\put(33.7000,-6.8000){\makebox(0,0)[lb]{1}}%
% STR 2 0 3 0 Black White  
% 4 2170 1930 2170 2030 2 0 0 0
% 1
\put(21.7000,-20.3000){\makebox(0,0)[lb]{1}}%
% STR 2 0 3 0 Black White  
% 4 2580 1940 2580 2040 2 0 0 0
% 2
\put(25.8000,-20.4000){\makebox(0,0)[lb]{2}}%
% STR 2 0 3 0 Black White  
% 4 2970 1940 2970 2040 2 0 0 0
% 1
\put(29.7000,-20.4000){\makebox(0,0)[lb]{1}}%
\end{picture}}%
\end{center}
\bigskip
$\mathcal{F}(\widehat{G} \cup \widehat{G_0}) 
= \{(Ch(\pi_i), \circ), (Ch(\tau_j), \bullet): \pi_i \in \widehat{G}, \tau_j \in \widehat{G_0}\}$. 
Put $\gamma_i = (Ch(\pi_i), \circ)$ and $\rho_j = (Ch(\tau_j), \bullet)$.  
Then the structure equations are
\begin{align*}
&\gamma_1  \gamma_1 = \gamma_0,~~~
\gamma_2  \gamma_2 = \gamma_0 + \gamma_1 + \gamma_2,~~~
\gamma_3 \gamma_3 = \gamma_4  \gamma_4 
= \gamma_0 + \gamma_2 + \gamma_3 + \gamma_4,\\
&\gamma_1  \gamma_2 = \gamma_2,~~~\gamma_1  \gamma_3 = \gamma_4,~~~
\gamma_1 \gamma_4 = \gamma_3,~~~
\gamma_2 \gamma_3 = \gamma_2  \gamma_4 
= \gamma_3 + \gamma_4,\\
&\gamma_3 \gamma_4 
= \gamma_1 + \gamma_2 + \gamma_3 + \gamma_4,~~~
\rho_0 \rho_0 = \rho_1 \rho_1 = \gamma_0 + \gamma_3,\\
& \rho_2 \rho_2 = \gamma_0 + \gamma_1 + \gamma_2 + 
\gamma_3 + \gamma_4, ~~~
\rho_1 \rho_2 = \gamma_2 + \gamma_3 + \gamma_4,\\
&\gamma_0 \rho_0 = \rho_0,~~~\gamma_1 \rho_0 = \rho_1,~~~\gamma_2 \rho_0 = \rho_2,~~~
\gamma_3 \rho_0 = \rho_0 + \rho_2,~~~
\gamma_4 \rho_0 = \rho_1 + \rho_2,\\
&\gamma_0 \rho_1 = \rho_1,~~~\gamma_1 \rho_1 = \rho_0,~~~\gamma_2 \rho_1 = \rho_2,~~~
\gamma_3 \rho_1 = \rho_1 + \rho_2,~~~
\gamma_4 \rho_1 = \rho_0 + \rho_2,\\
&\gamma_0 \rho_2 = \rho_2,~~~\gamma_1 \rho_2 = \rho_2,~~~
\gamma_2 \rho_2 = \rho_0 + \rho_1 + \rho_2,\\
&\gamma_3 \rho_2 = \gamma_4 \rho_0 = \rho_0 + \rho_1 + \rho_2. 
\end{align*}
As in the case that example 4.1, we obtain a hypergroup isomorphic to 
the hypergroup of [HKTY1, example 4.10]. We write that only dimension function of 
$\mathcal{F}(\widehat{G} \cup \widehat{G_0})$. 
\begin{align*}
&d(\gamma_0) = d(\gamma_1) = 1,~~~d(\gamma_2) = 2,~~~d(\gamma_3) = d(\gamma_4) = 3,\\
&d(\rho_0) = d(\rho_1) = 2,~~~d(\rho_2)= 4. 
\end{align*}

\bigskip
\noindent
{\bf Acknowledgment } The authors would like to thank Satoshi Kawakami 
who gave us many valuable comments and warm encouragement.

\bigskip

\textbf{Addresses }

\medskip

Narufumi Nakagaki : Nara University of Education

Graduate School of Education

Takabatake-cho

Nara, 630-8528

Japan

\medskip

%e-mail : a153304@student.nara-edu.ac.jp

\bigskip

Tatsuya Tsurii : Osaka Prefecture University 

1-1 Gakuen-cho, Nakaku, Sakai 

Osaka, 599-8531

Japan

\medskip

e-mail : tatsuya\_t\_99@yahoo.co.jp

\end {document}